\documentclass[]{elsart}
\usepackage{amsmath,amssymb%
}
\def\bbE{{\mathbb E}}

\def\bbN{{\mathbb N}}

\def\R{\Re}

\def\bP{{\bf P}}
\def\bE{{\bf E}}

\def\bT{\mbox{\bf T}}
\def\bQ{\mbox{\bf Q}}
\def\ovpi{\overline{\pi}}
\def\ovbet{\overline{b}}
\def\ovgamma{\overline{\gamma}}
\def\ovdelta{\overline{\delta}}
\def\ovp{\overline{p}}
\def\ovpi{\overline{\pi}}
\def\ovb{\overline{b}}
\def\ovf0{\widehat{f}^0}
\def\ubeta{\underline{\beta}}
\def\ux{\underline{x}}
\def\uy{\underline{y}}

\def\uX{\underline{X}}
\def\ovPi{\overline{\Pi}}
\def\ovBet{\overline{B}}

\def\itPsi{\mathit{\Psi}}
\def\itPsiT{\mathit{\widetilde{\Psi}}}
\def\ST{\widetilde{S}}
\def\itLambda{\mathit{\Lambda}}
\def\itLambdaT{\mathit{\widetilde{\Lambda}}}
\def\itXi{\mathit{\Xi}}

\def\phiP{{\bP^\varphi}}
\def\phiE{\bE^\varphi}

\def\cF{{\mathcal F}}

\def\cB{{\mathcal B}}

\def\u{{\sf u}}

\def\frF{{\mathfrak F}}
\def\frS{{\mathfrak S}}

\def\hwtilde{\widetilde h}

\def\one{{\mathbb I}}

\newtheorem{ptheorem}{Theorem}
\newtheorem{lemma}{Lemma}

\newtheorem{premark}{Remark}

\usepackage[cp1250]{inputenc}
\usepackage{fancyhdr}
\usepackage{color}
\usepackage[]{hyperref}

\begin{document}
\small
\begin{frontmatter}

\title{Unspecified distribution in single disorder problem}%\thanksref{label1}}
\thanks[label1]{\today}

%\thanks[lab1]{The research was supported by KBN grant no  2 P03A 021 22 (350228).}
\author{Wojciech Sarnowski\thanksref{label2}\thanksref{cor1}}
\ead{Wojciech.Sarnowski@pwr.wroc.pl}
\ead[url]{http://www.im.pwr.wroc.pl/\~{}sarnowski}
\corauth[cor1]{Corresponding author}
\author{Krzysztof Szajowski\thanksref{label2}}%\corauthref{cor1}}
\ead{Krzysztof.Szajowski@pwr.wroc.pl}
\ead[url]{http://neyman.im.pwr.wroc.pl/\~{}szajow}
\address[label2]{Wroc\l{}aw University of Technology, Institute of
Mathematics and Computer Science, Wybrze\.{z}e Wyspia\'{n}skiego 27, 50-370 Wroc\l{}aw, Poland}

\date{ \today }
\maketitle

\begin{abstract}
We register a stochastic sequence affected by one disorder. Monitoring of the sequence is made in the
circumstances when not full information about distributions before and after the change is available. The initial problem
of disorder detection is transformed to optimal stopping of observed sequence. Formula for optimal decision functions
is derived.
\end{abstract}

%\PACS
%\MSC Primar 60G40 \sep 60K99; \quad Secondary 90D60

\end{frontmatter}
%%%%%%%%%%%%%%%%%%%%%%%%%%%%%%%%%%%%%%%%%%%%%%%%%%%%%%%%%%%%%%%%%%%%%%%%%%%%%%%%%%%%%%%%%%%%%%%%%%%
\fancyhf{}
\fancyhead[FC]{\thepage}
\fancyhead[HC]{}
\thispagestyle{fancy}
%   { \LARGE {\bf
%    Probability maximizing approach to double disorder problem}}

\section{Introduction}
The paper is focused on sequential detection using Bayesian approach. Disorder problem in this framework was formulated
by A.N. Kolmogorov at the end of 50's of previous century and solved by ~\cite{shi61:detection}. The next turning point
is paper by ~\cite{pesshi02:poisson} where authors provide complete solution of basic problem. From this time many
publications provide new solutions and generalizations in the area of sequential detection. Some of them are articles
by ~\cite{kar03:note} and ~\cite{baydaykar05:poisson}. For discrete time case there are some detailed analysis in the
papers by %~\cite{boj79:disorder},
~\cite{bojhos84:problem}, ~\cite{mou98:abrupt}, ~\cite{yak94:finite}, ~\cite{yos83:complicated}, ~\cite{sza96:twodis} and the papers cited there.

Such model of data appears in many practical problems of the quality control~(see Brodsky and Darkhovsky~\cite{brodar93:nonparametr}, Shewhart~\cite{she31:quality} and in the collection of the papers \cite{basben86:abrupt}), traffic anomalies in networks (in papers by Dube and Mazumdar~\cite{dubmaz01:quickest}, Tartakovsky et al.~\cite{tarroz06:intrusions}), epidemiology models (see Baron~\cite{bar04:epidemio}). In management of manufacture it happens that the plants which produce some details changes their parameters. It makes that the details change their quality. The aim is to recognize the moments of these changes as soon as possible.

This paper focuses attention on models under assumption of uncertainty about distribution
before or after the change. The example of such models can be found in research by ~\cite{dubmaz01:quickest} with
application to detection of traffic anomalies in networks or in paper by ~\cite{sarsza08:disorder}. The solution of a
single disorder model with unspecified distribution of observed sequence is presented. Section~\ref{PojRozregNieznRozkl-Model} specifies the details of investigated model. The transformation of the optimization job to the optimal stopping problem for the specific stochastic process is considered in Section~\ref{PojRozregNieznRozkl-IstnienieRozwiazania}. A construction of the optimal estimator of the disorder moment is given in Section~\ref{PojRozregNieznRozkl-Oznaczenia}. Technical parts of investigations are moved to Appendix.

\section{Description of the model}\label{PojRozregNieznRozkl-Model}%Section 2
\subsection{\label{PojRozregNieznRozkl-OznaczeniaWektory}Basic notations}

For further considerations it will be convenient to introduce the following notation which will make our formulas
more compact and clear
\begin{eqnarray}
    \underline{x}_{k,n} &=& (x_k, x_{k+1},...,x_{n-1},x_n),\; k\leq n,\nonumber\\
     L_{m}^{i,j}(\underline{x}_{k,n}) &=& \prod_{r=k+1}^{n-m} \!\!f_{x_{r-1}}^{0,i}(x_{r})\!\!\!\! \prod_{r=n-m+1}^{n} \!\!\!f_{x_{r-1}}^{1,j}(x_{r}), \nonumber\\ %%%; \;\;k \leq m \leq n-k. \nonumber
    \underline{A}_{k,n} &=& A_k \times A_{k+1} \times \ldots \times A_n, \nonumber
\end{eqnarray}
where: $\prod_{r=m_1}^{m_2}u_r = 1$ for $m_1 > m_2$ and $u_r \in \R$, $A_i \in \cB$, $k \leq i \leq n$.

It will be convenient to write $\ubeta = (\beta_1,\beta_2)$ and denote by $\overline{\alpha} = \left(\alpha_{11},\ldots,\alpha_{1l_1},\ldots, \alpha_{l_01},\ldots,\alpha_{l_0l_1} \right)$
any matrix $l_0 \times l_1$:
\begin{eqnarray}
\left[
            \begin{array}{cccc}
                \alpha_{11} & \alpha_{12} & \cdots & \alpha_{1l_1} \\
                \alpha_{21} & \alpha_{22} & \cdots & \alpha_{2l_1} \\
                \vdots & \vdots & \ddots & \vdots  \\
                \alpha_{l_01} &  \alpha_{l_0 2} & \cdots & \alpha_{l_0l_1} \\
            \end{array}
\right]\nonumber
\end{eqnarray}
In consequence vectors $\ovpi$, $\ovbet$, $\ovp$ represent:
\begin{eqnarray}
\ovpi &=& \left(\pi_{11},\ldots,\pi_{1l_2},\ldots, \pi_{l_01},\ldots,\pi_{l_0l_1} \right) \nonumber\\
\ovbet &=& \left(b_{11},\ldots,b_{1l_2},\ldots, b_{l_01},\ldots,b_{l_0l_1} \right) \nonumber\\
\ovp &=& \left(p_{11},\ldots,p_{1l_1},\ldots, p_{l_01},\ldots,p_{l_0l_1} \right) \nonumber
\end{eqnarray}

We need also notation for vector of densities $f_x^{0,i}(y)$. Let $\ovf0_x(y)$, where $x,y \in \bbE$ stands
behind:
\begin{eqnarray}
    \ovf0_x(y) = ( \underbrace{f_x^{0,1}(y),\ldots,f_x^{0,1}(y) }_{ l_1\; \rm{times}},\ldots,  \underbrace{f_x^{0,l_0}(y),\ldots,f_x^{0,l_0}(y) }_{ l_1\; \rm{times}} ). \nonumber
\end{eqnarray}
Moreover let us introduce operation "$\circ$". For vectors $\overline{\alpha}$ and $\overline{\beta}$ we put:
\begin{eqnarray}
\overline{\alpha}\circ \overline{\beta}= \left(\alpha_{11}\beta_{11},\ldots,\alpha_{1l_1}\beta_{1l_1},\ldots, \alpha_{l_01}\beta_{l_01},\ldots,\alpha_{l_0l_1}\beta_{l_0l_1} \right).\nonumber
\end{eqnarray}

\subsection{Change point problem}
Let $(X_n)_{ n \in \bbN}$ be sequence of observable random variables defined on $(\Omega,\cF, \bP)$ with value in $(\bbE, \cB)$, $\bbE \subset\Re$. Sequence $(X_n)$ generates filtration $\cF_n = \sigma(X_0,X_1,...,X_n)$. On the same space there are also defined  variables $\theta$, $\beta_1$ and $\beta_2$. $\theta$ takes values in $\{1,2,3,\ldots\}$. Variables $\beta_1$, $\beta_2$ are valued in $I_k = \{1,2,\ldots, l_k\}$, where $l_k \in~\bbN$, $k=0,1$. Let us assume the following parametrization:
\begin{eqnarray*}
\label{PojRozregNieznRozkl-rozkladBeta}
\bP(\beta_1=i, \beta_2 = j) &=& b_{ij} \\
\label{PojRozregNieznRozkl-rozkladTetaWarBeta}
\bP(\theta = n | \beta_1=i, \beta_2 = j) &=& \left\{
    \begin{array}{ll}
        \pi_{ij},                       & \mbox{if $n=1$,}\\
        (1-\pi_{ij})p_{ij}^{n-2}q_{ij}, & \mbox{if $n > 1$,}\\
    \end{array}
\right.
\end{eqnarray*}
where $i \in I_0, j\in I_1$, $\sum_{i \in I_0, j\in I_1}b_{ij}=1$, $b_{ij} \geq 0$, $\pi_{ij} \in [0,1]$,
$p_{ij} = 1 - q_{ij} \in (0,1)$.
%{
%\setlength\arraycolsep{0.2pt}
We have
\begin{equation*}
    \sum_{k=1}^{\infty} \sum_{i \in I_0} \sum_{j \in I_1}\bP(\theta = k, \beta_1=i, \beta_2 = j)=1
\end{equation*}
The change of the conditional densities in random moment $\theta$ is investigate in this model. The transfer between
distribution is described by conditional probabilities $b_{ij}=\bP(\beta_2=j|\beta_1=i)$. For completeness it will be
assumed that the state of $\beta_1$ is stable before $\theta$ and the same as at the moment $0$. The marginal distribution
of $\theta$ has a form
\begin{eqnarray*}
\label{PojRozregNieznRozkl-rozkladBrzegTeta}
\bP(\theta = k) &=& \sum_{i,j}\bP(\theta = k, \beta_1=i, \beta_2 = j) \nonumber \\
&=& \left\{
    \begin{array}{ll}
        \sum_{i,j}\pi_{ij}\cdot b_{ij}             & \mbox{if $k=1$,}\\
        \sum_{i,j}(1-\pi_{ij})p_{ij}^{k-2}q_{ij}b_{ij} & \mbox{if $k > 1$.}\\
    \end{array}
\right.
\end{eqnarray*}
The observed sequence has a form
\begin{eqnarray}
\label{PojRozregNieznRozkl-procesyX}
    X_n = X^{0,i}_n \cdot \one_{\{\theta > n,\; \beta_1 = i\}} + X^{1,j}_n \cdot \one_{\{\theta \leq n,\; \beta_2 = j\;, X^{1,j}_{\theta-1}=X^{0,i}_{\theta-1} \}},
\end{eqnarray}
where $(X_n^{r,i}, \mathcal{G}_n^{r,i}, \bP_x^{r,i})$, $r=0,1$, are Markov processes and $\sigma$-fields: $\mathcal{G}_n^{r,i} = \sigma(X_0^{r,i}, X_1^{r,i},\ldots ,X_n^{r,i})$, with $i \in I_0$, $j \in I_1$, $r=0,1$ and
$n \in \{ 0, 1, 2, \ldots\}$. Variables $\theta$, $\beta_1$ and $\beta_2$ are not measurable w.r.t $\cF_n$.

On the space $(\bbE, \cB)$ there are $\sigma$-additive measures $\mu(\cdot)$ and measures $\mu_x^{\bullet,\bullet}$ absolutely continuous with respect to $\mu$. It is assumed that the measures $\bP_x^{k,i}(\cdot)$, $i=1,2,\ldots,l_k$, $k=0,1$, have following representation:
\begin{eqnarray*}
\bP_x^{k,i}(\{\omega:X_1^{k,i}\in B\})&=&\bP(X_1^{k,i}\in B|X_0^{k,i}=x)=\int_Bf_x^{k,i}(y)\mu(dy)\\
&=&\int_B\mu_x^{k,i}(dy)=\mu_x^{k,i}(B).
\end{eqnarray*}
for any $B\in\cB$. The conditional densities $f_x^{k,1}(\cdot),\ldots, f_x^{k,l_k}(\cdot)$ are different and supports of all measures $\mu_x^{\cdot,\cdot}$ are there same for given $x\in\bbE$.
It is the model of the following random phenomenon. At the beginning we register process $\{X_n^{0,i},n \in \mathbb{N} \}$, where $i\in I_0$ is unknown. At random moment $\theta$ initial process is switched on $\{X_n^{1,j},n \in \mathbb{N} \}$ where $j\in I_1$ is unknown. It can be interpreted as disorder of $\{X_n,n \in \bbN \}$ causing change in distribution of $\{X_n\}_{n \in \bbN}$. We monitor the process and we wish to detect the change as close $\theta$ as possible. However our knowledge about densities before and after the change moment $\theta$ is limited generally to the information about sets of possible conditional densities only:
$ \{f_x^{0,i}(y), i \in I_0 \}$ and $\{f_x^{1,j}(y), j \in I_1\}$
respectively. We also know probabilities of distribution pairs $b_{ij}$ and parameters $\pi_{ij}$.

For $i \in I_0, j \in I_1$ let us introduce functions ${\itPsi}^{i,j}, {\itPsiT}^{i,j}, \itLambda^{i,j}, \itLambdaT^{i,j}$
defined on the product $\bbN \times (\times_{i=1}^{l+2} \bbE) \times [0,1]$ with values in $\Re$:
{
\small
\begin{eqnarray}\label{denomAuxFunctPi}
\!\!\!\!\!\!\!\!\!\!\!\!\!\!\!\!\!\!\!\!\!\!\!\!\!{\itPsi}^{i,j}(l,\ux_{0,l+1},\alpha)&=& (1-\alpha)\!\!\left[ q_{ij}\!\!\sum_{k=0}^l p_{ij}^{l-k}L_{k+1}^{i,j}(\underline{x}_{0,l+1}) + p_{ij}^{l+1}L_{0}^{i,j}(\underline{x}_{0,l+1}) \right]+\alpha L_{l+1}^{i,j}(\underline{x}_{0,l+1})\\
%\nonumber&&\mbox{}+\alpha L_{l+1}^{i,j}(\underline{x}_{0,l+1}), \\
\label{denomAuxFunctPiTilda}
\!\!\!\!\!\!\!\!\!\!\!\!\!\!\!\!\!\!\!\!\!\!\!\!\!{\itPsiT}^{i,j}(l,\ux_{0,l+1},\alpha)&=& (1-\alpha)\!\!\left[ q_{ij}\!\!\sum_{k=1}^l p_{ij}^{l-k}L_{k}^{i,j}(\underline{x}_{0,l+1}) + p_{ij}^{l}L_{0}^{i,j}(\underline{x}_{0,l+1})\right]+\alpha L_{l+1}^{i,j}(\underline{x}_{0,l+1}),  \\
%\nonumber&&\mbox{}+\alpha L_{l+1}^{i,j}(\underline{x}_{0,l+1}),\\
\nonumber%\label{auxFunctLambda}
\!\!\!\!\!\!\!\!\!\!\!\!\!\!\!\!\!\!\!\!\!\!\!\!\!\itLambda^{i,j}(l,\ux_{0,l+1},\alpha) &=&  {\itPsi}^{i,j}(l,\ux_{0,l+1},\alpha) - (1-\alpha)p_{ij}^{l+1}L_{0}^{i,j}(\underline{x}_{0,l+1}),\\
\nonumber%\label{auxFunctLambdaTilda}
\!\!\!\!\!\!\!\!\!\!\!\!\!\!\!\!\!\!\!\!\!\!\!\!\!\itLambdaT^{i,j}(l,\ux_{0,l+1},\alpha) &=& {\itPsiT}^{i,j}(l,\ux_{0,l+1},\alpha) - (1-\alpha)p_{ij}^{l}L_{0}^{i,j}(\underline{x}_{0,l+1}).
\end{eqnarray}
}
\normalsize
Next let us define on $\bbN \times (\times_{i=1}^{k+2} \bbE) \times (\times_{i=1}^{l_1l_2}[0,1]) \times (\times_{i=1}^{l_1l_2}[0,1])$ function $S, \ST$:
\begin{eqnarray}
\label{multidist}
S(k,\ux_{0,k+1}, \overline{\gamma}, \overline{\delta})&=&\sum_{i,j}\gamma_{ij}\itPsi^{i,j}(k,\ux_{0,k+1},\delta_{ij}),\\
\label{multidistTilda}
\ST(k,\ux_{0,k+1}, \overline{\gamma}, \overline{\delta})&=&\sum_{i,j}\gamma_{ij}\itPsiT^{i,j}(k,\ux_{0,k+1},\delta_{ij}).
\end{eqnarray}

For any $D_n=\{\omega:X_i\in B_i$, $i=1,2,\ldots,n\}$, where $B_i\in\cB$ and any $x\in\bbE$ define:
\begin{eqnarray*}
\bP_x(D_n)=\bP(D_n|X_0=x)&=&\int_{\times_{i=1}^nB_i}\ST(n-1,\ux_{0,n},\ovb,\ovpi)\mu(d\ux_{1,n})%\\
%&=&\int_{\times_{i=1}^nB_i}\mu_{x_0}(d\ux_{1,n})=\mu_{x_0}(\times_{i=1}^nB_i).
\end{eqnarray*}

For the process (\ref{PojRozregNieznRozkl-procesyX}) the set of estimators for the disorder moment $\theta$ is  $\frS^X$ -- the set of stopping times with respect to $\{\cF_n\}_{n \in \bbN\cup\{0\}}$. The construction of the optimal estimator is to find a stopping time $\tau^{*}\in \frS^X$ such that for any $x\in\bbE$
\begin{equation}
\label{PojRozregNieznRozkl-Problem}
  \bP_x( | \theta - \tau^{*} | \leq d ) = \sup_{\tau \in \frS^X} \bP_x( | \theta - \tau | \leq d ),
\end{equation}
where $d \in \{0,1,2,...\}$ is fixed level of detection precision.

\section{Existence of solution}\label{PojRozregNieznRozkl-IstnienieRozwiazania}%Section 3
In this section we are going to show that there exists solution of the problem (\ref{PojRozregNieznRozkl-Problem}).
Let us define:
\begin{eqnarray}
Z_n &=& \bP(| \theta - n | \leq d \mid \cF_n),\; n=1,2,\ldots, \nonumber \\
V_n &=& \rm{ess}\sup_{\{\tau \in \frS^X,\;\tau \geq n\}}\bP( | \theta - n | \leq d \mid \cF_n), \; n=0,1,2,\ldots \nonumber \\
\label{PojRozregNieznRozkl-StopIntuicyjny}
    \tau_0 &=& \inf\{ n: Z_n=V_n \}
\end{eqnarray}
Notice that, if $Z_{\infty}=0$, then $Z_{\tau} = \bP( |\theta - \tau | \leq d \mid \cF_{\tau})$for
$\tau \in \frS^X$. Because $\cF_{n} \subseteq \cF_{\tau}$ (when $n \leq \tau$), we obtain
\begin{eqnarray}
    V_n &=& \rm{ess}\sup_{\tau \geq n}\bP(|\theta - \tau | \leq d \mid \cF_n) = \rm{ess}\sup_{\tau \geq n} \bE(\bE(\one_{\{ |\theta - \tau | \leq d\}}|\cF_\tau) \mid \cF_n) \nonumber \\
            &=& \rm{ess}\sup_{\tau \geq n}\bE(Z_{\tau} \mid \cF_n) \nonumber
\end{eqnarray}
The following lemma states that solution exists.
\begin{lemma}
    \label{PojRozregNieznRozkl-lematCzasStopu}
    Stopping time $\tau_0$ given by (\ref{PojRozregNieznRozkl-StopIntuicyjny}) is a solution of the problem
    (\ref{PojRozregNieznRozkl-Problem}).
\end{lemma}
\begin{pf} Applying Theorem~1 from \cite{boj79:disorder} it is enough to show that
$\displaystyle{\lim_{n \rightarrow \infty}Z_n=0}$. For all $n,k$, where $n\geq k$ we have:
\begin{eqnarray*}
   Z_n &=& \bE(\one_{\{ |\theta - n | \leq d \}} \mid \cF_n) \leq \bE(\sup_{j \geq k}\one_{\{ |\theta - j | \leq d \}} \mid \cF_n)
\end{eqnarray*}
%\normalsize
Basing on Levy's theorem we get
$\limsup_{n\rightarrow \infty}Z_n \leq \bE(\sup_{j \geq k}\one_{\{ |\theta - j| \leq d \}} \mid \cF_{\infty})$
where $\cF_{\infty} = \sigma\left( \bigcup_{n=1}^{\infty}\cF_n \right)$.
We have: $\limsup_{j \geq k,\; k\rightarrow \infty}\one_{\{ |\theta - j | \leq d \}} = 0$ \emph{a.s.}
Basing on dominated convergence theorem we get we state that
\[
  \lim_{k \rightarrow \infty}\bE(\sup_{j\geq k}\one_{\{ |\theta - j | \leq d \}} \mid \cF_{\infty} ) = 0\;\; a.s.
\]
what ends the proof.
\end{pf}
It turns out that we need at least $d$ observations to detect disorder in optimal way:
\begin{lemma}
    \label{PojRozregNieznRozkl-lematCzasStopuMax}
Let $\tau$ be stopping rule in the problem (\ref{PojRozregNieznRozkl-Problem}). Then rule
$\tilde{\tau} = \max(\tau,d+1)$ is at least as good as $\tau$ (in the sense of (\ref{PojRozregNieznRozkl-Problem})).
\end{lemma}
\begin{pf} For $\tau \geq d+1$ the rules are the same. Let us consider case when $\tau < d+1$. Then $\tilde{\tau} = d+1$ and:
\begin{eqnarray}
\bP(|\theta - \tau| \leq d) &=& \bP(\tau - d  \leq \theta \leq \tau + d)
                          = \bP( 1 \leq \theta \leq \tau + d)                         \nonumber \\
                          &\leq& \bP( 1 \leq \theta \leq 2d+1)
                          = \bP(\tilde{\tau} - d  \leq \theta \leq \tilde{\tau} + d)
                          = \bP(|\theta - \tilde{\tau}| \leq d).                      \nonumber
\end{eqnarray}
\end{pf}

\section{Construction of the disorder moment estimator}\label{PojRozregNieznRozkl-Oznaczenia}%Section 4
\subsection{Function and processes}%\label{PojRozregNieznRozkl-OznaczeniaFunProc}
Let us fix parameters $\ovpi$, $\ovbet$ and set initial state of $X_n$: $\bP(X_0 = x ) = 1$.
%\begin{eqnarray}
%\bP^{\ovpi,\ovbet}(X_0 = x ) = 1 \nonumber
%\end{eqnarray}
We denote $\varphi = (\ovpi,\ovbet,x)$ and we will write $\bP^\varphi(\bullet)$ to emphasis that the probability of the
events defined by the process are dependent on this \emph{a priori} set parameters. Let us define the following crucial
posterior processes:
\begin{eqnarray}
\label{ProcesPin}
\Pi_n^{i,j} &=& \phiP(\theta \leq n | \ubeta = (i,j), \cF_n )
             =  \phiP(\theta \leq n | \tilde{\cF}_n^{ij} ) \\
\label{ProcesBetan}
B_n^{i,j} &=& \bP^{\varphi}(\ubeta = (i,j) | \cF_n )
\end{eqnarray}
where $n \in \bbN$, $i \in I_0, j \in I_1$, $\tilde{\cF}_n^{i,j}=\sigma(\cF_n,\one_{\{\ubeta=(i,j)\}})$. Process $\Pi_n^{i,j}$ is designed for updating information about disorder
distribution. $B_n^{i,j}$ in turn refreshes information about distributions of variables $\beta_1$, $\beta_2$.
Notice that $\Pi_n^{i,j}$, $B_n^{i,j}$ starts from following states: $\Pi_0^{i,j} = 0$, $B_0^{i,j} = b_{ij}$
%(for proof see section \ref{PojRozregNieznRozkl-dodatekWzor4} ).
Dynamics of $\Pi_n^{i,j}$ and $B_n^{i,j}$ are characterized by formulas (\ref{PojRozregNieznRozkl-Pi_n_Rekurencja2}),
(\ref{PojRozregNieznRozkl-B_n_Rekurencja}). The above notations hold also for (\ref{ProcesPin}), (\ref{ProcesBetan}):
\begin{eqnarray}
\ovPi_n &=& \left(\Pi_n^{1,1},\ldots,\Pi_n^{1,l_2},\ldots, \Pi_n^{l_1,1},\ldots,\Pi_n^{l_1l_2} \right), \nonumber\\
\ovBet_n &=& \left(B_n^{1,1},\ldots,B_n^{1,l_2},\ldots, B_n^{l_1,1},\ldots,B_n^{l_1l_2} \right). \nonumber
\end{eqnarray}
At the end of section let us define auxiliary functions $\Pi^{\cdot,\cdot}(\cdot,\cdot,\cdot)$,
$\Gamma^{\cdot,\cdot}(\cdot,\cdot,\cdot,\cdot)$. For $x, y \in \bbE$,
$\alpha, \gamma_{ij}, \delta_{ij} \in [0,1]$, $i \in I_0$, $j \in I_1$ put:
%\begin{eqnarray}
%\label{PojRozregNieznRozkl_funkcjaH3}
%    H^{i,j}(x,y,\alpha) = f_{x}^{1,j}(y)(q_{ij}+p_{ij}\alpha) + f_{x}^{0,i}(y)p_{ij}(1-\alpha).
%\end{eqnarray}
%and
\begin{eqnarray}\label{auxFunctPi}
\Pi^{i,j}(k,\ux_{0,n},\alpha) &=& \frac{\itLambda^{i,j}(k,\ux_{0,n},\alpha)}{{\itPsi}^{i,j}(k,\ux_{0,n},\alpha)}\\
\label{auxFunctGamma}
\Gamma^{i,j}(k,\ux_{0,n},\ovgamma, \ovdelta) &=&\frac{\gamma_{ij}{\itPsi}^{i,j}(k,\ux_{0,n},\delta_{ij})}{S(k,\ux_{0,n}, \ovgamma,\ovdelta)}.
\end{eqnarray}

Let $D_n=\{\omega:\underline{X}_{0,n}\in \underline{B}_{0,n}\}$, $X_0=x$ and $B_i\in\cB$. We have
\begin{eqnarray}\label{SnijFormula1}
\phiP(\theta>n,\beta=(i,j),D_n)&=&\int\limits_{\{\omega:\beta=(i,j),D_n\}}\one_{\{\theta>n\}}d\phiP\\
\nonumber&\hspace{-20em}=&\hspace{-10em}  \int\limits_{\underline{B}_{0,n}}\frac{(1-\pi_{ij})p_{ij}^{n-1}L_0^{ij}(\ux_{0,n})}{S_n^{i,j}(\ux_{0,n})}\frac{b_{ij}S_n^{i,j}(\ux_{0,n})}{S_n(\ux_{0,n})}S_n(\ux_{0,n})\mu(d\ux_{1,n})\\
\nonumber&\hspace{-20em}=&\hspace{-10em}\int\limits_{D_n}(1-\Pi_n^{i,j})B_n^{i,j}d\phiP,
\end{eqnarray}
where
\begin{eqnarray}
\label{Snij}
S_n^{i,j}(\ux_{0,n})&=&\pi_{ij}L_n^{i,j}(\ux_{0,n})+(1-\pi_{ij})p_{ij}^{n-1}L_0^{i,j}(\ux_{0,n})\\
\nonumber&&\mbox{}+(1-\pi_{ij})\sum_{s=2}^np_{ij}^{s-2}q_{ij}L_{n-s+1}^{i,j}(\ux_{0,n})=\Psi^{i,j}(n-1,\ux_{0,n},\pi_{ij})
\end{eqnarray}
and $S_n(\ux_{0,n})=\sum_{i,j}b_{ij}S_n^{i,j}(\ux_{0,n})=S(n-1,\ux_{0,n},\bar{b},\bar{\pi})$.
%\begin{equation}\label{auxFunctLambda}
%\itLambda^{i,j}(d,\ux_{0,d+1},\alpha)=\alpha L_{d+1}^{i,j}(\underline{x}_{0,d+1})+(1-\alpha)q_{ij}\sum_{k=0}^d p_{ij}^{d-k}L_{k+1}^{i,j}(\underline{x}_{0,d+1}).
%\end{equation}
%It will be convenient to have also for $n > s+1$
%\begin{eqnarray}\label{auxFunctUpsilon}
%\itUpsilon^{i,j}(s,\ux_{n-s-1,n},\alpha)&=&\alpha L_{s+1}^{i,j}(\underline{x}_{n-s-1,n})\\
%&&+(1-\alpha)\big[q_{ij}\sum_{k=0}^s p_{ij}^{s-k}L_{k+1}^{i,j}(\underline{x}_{n-s-1,n})+p_{ij}^{s+1} L_{0}^{i,j}(\underline{x}_{n-s-1,n})\big].
%\end{eqnarray}

%%%%%%%%%% Solution of auxiliary optimal stopping problem %%%%%%%%%%%%%%%%%%%%%%%%%%%%%%%%%%%%%%%%%%%
\section{Solution}\label{PojRozregNieznRozkl-RozwMetodShiryaev}
According to Shiryayev's methodology (see \cite{shi78:optimal} ) we are going to find solution reducing initial
problem (\ref{PojRozregNieznRozkl-Problem}) to the case of stopping Random Markov Function with special payoff
function. This will be done using posterior processes (\ref{ProcesPin})-(\ref{ProcesBetan}).
\begin{lemma}
    \label{PojRozregNieznRozkl-ZmianaWyplaty}
For $n\geq d+1$
\begin{eqnarray}
%\label{PojRozregNieznRozkl-PodstWyplataJakoWWO}
\label{PojRozregNieznRozkl-StaraWyplataNowaWyplata1}
\phiP(| \theta - n | \leq d) =
\left\{
    \begin{array}{ll}
      \phiE\left[ h(\underline{X}_{n-1-d,n},\ovPi_n,\ovBet_n) \right], & \mbox{if $n >d+1$},\\
        &\\
      \phiE\left[ \widetilde{h}(\ovPi_{d+1},\ovBet_{d+1}) \right], &  \mbox{if $n = d+1$}.\\
    \end{array}
            \right.
\end{eqnarray}
where
\begin{eqnarray}
\label{PojRozregNieznRozkl-NowaWyplata}
h({\ux}_{1,d+2},\overline{\gamma}, \overline{\delta}) &=& \sum_{i,j}\left( 1-p_{ij}^d+ q_{ij}\sum_{k=1}^{d+1}\frac{L^{i,j}_{k}(\underline{x}_{1,d+2})} { p_{ij}^k\ L^{i,j}_{0}(\underline{x}_{1,d+2})} \right) (1-\gamma_{ij})\delta_{ij},\\
\label{PojRozregNieznRozkl-NowaWyplata_dplus1}
\hwtilde(\overline{\gamma}, \overline{\delta}) &=& \sum_{i,j}\left( 1-p_{ij}^d(1-\gamma_{ij})\right)\delta_{ij},
\end{eqnarray}
$x_1,...,x_{d+2} \in \bbE$, $\gamma_{ij},\delta_{ij} \in [0,1]$, $i \in I_0$, $j \in I_1$.
\end{lemma}
\begin{pf} Let us rewrite initial criterion as expectation:
\begin{eqnarray}
\label{PojRozregNieznRozkl-PodstWyplataJakoWWO}
\phiP(| \theta - n | \leq d) &=& \phiE\left[ \phiP(| \theta - n | \leq d \mid \cF_n ) \right].
\end{eqnarray}
Let us analyze conditional probability under expectation in equation (\ref{PojRozregNieznRozkl-PodstWyplataJakoWWO})
using total probability formula
{%\setlength\arraycolsep{0pt}
\begin{eqnarray}
\label{PojRozregNieznRozkl-PodstWyplataJakoWWO1}
\!\!\!\!\!\!\!\phiP&(&| \theta - n | \leq d \mid \cF_n ) =  \phiP( \theta \leq n+d \mid \cF_n ) - \phiP( \theta \leq n-d-1 \mid \cF_n ) \\
\!\!\!\!\!\!\!&=&  \sum_{i,j}\Pi^{ij}_{n+d}B_n^{i,j}- \sum_{i,j}\Pi^{ij}_{n-d-1}B_n^{i,j},\nonumber
\end{eqnarray}}
because
\begin{eqnarray*}
\phiP(\theta\leq n+d|\cF_n)&=&\phiE(\one_{\theta\leq n+d}|\cF_n)=\sum_{i,j}\phiE(\one_{\{\theta\leq n+d\}}\one_{\{\ubeta=(i,j)\}}|\cF_n)\\
&=&\sum_{i,j}\phiE(\phiE(\one_{\{\theta\leq n+d\}}\one_{\{\ubeta=(i,j)\}}|\tilde{\cF}_n^{i,j})|\cF_n)\\
&=&\sum_{i,j}\phiE(\one_{\{\ubeta=(i,j)\}}\phiE(\one_{\{\theta\leq n+d\}}|\tilde{\cF}_n^{i,j})|\cF_n)\\
&=&\sum_{i,j}\phiE(\one_{\{\theta\leq n+d\}}|\tilde{\cF}_n^{i,j})\phiE(\one_{\{\ubeta=(i,j)\}}|\cF_n).
\end{eqnarray*}
The last equality is a consequence of the very special form of the extended $\sigma$-field $\tilde{\cF}_n^{i,j}$. The
random variable measurable with respect to $\tilde{\cF}_n^{i,j}$ is also $\cF_n$ measurable. Putting $n=d+1$ in Lemma \ref{PojRozregNieznRozkl-dodatekWzor3} we get $\phiP( \theta \leq n-d-1 \mid \cF_n, \ubeta = (i,j) )=0$, for $i \in I_0, j \in I_1$. Hence
\[
\phiP(| \theta - n | \leq d \mid \cF_n ) = \sum_{i,j}\phiP( \theta \leq n+d \mid \tilde{\cF}_n )\phiP(\ubeta = (i,j)\mid {\cF}_n).
\]
Lemma \ref{PojRozregNieznRozkl-dodatekWzor2} implies that
$    \phiP(| \theta - n | \leq d) = \phiE\left[ \widetilde{h}(\ovPi_{d+1},\ovBet_{d+1})\right]$.
Now let $n>d+1$. Basing on Lemma \ref{PojRozregNieznRozkl-dodatekWzor2} probability
$\phiP( \theta \leq n+d \mid \tilde{\cF}_n )$ is given by (\ref{PojRozregNieznRozkl-PierwszyCzlonFunkcjiWyplaty}).
From Lemma \ref{PojRozregNieznRozkl-dodatekWzor3} we know that $\phiP( \theta \leq n-d-1 \mid \tilde{\cF}_n )$
is expressed by equation (\ref{PojRozregNieznRozkl-DrugiCzlonFunkcjiWyplaty1}). Formula
(\ref{PojRozregNieznRozkl-DrugiCzlonFunkcjiWyplaty1}) reveals connection between payoff function
(\ref{PojRozregNieznRozkl-StaraWyplataNowaWyplata1}) and posterior process at instants $n$ and $n-d-1$, i.e $\Pi_n^{i,j}$, $\Pi_{n-d-1}^{i,j}$ for $i \in I_0, j \in I_1$. Dependence on $\Pi_{n-d-1}^{i,j}$ can be rule out by expressing $\Pi_{n-d-1}^{i,j}$ in terms of $\Pi_{n}^{i,j}$. By Lemma~\ref{PojRozregNieznRozkl-dodatekWzor4} and (\ref{PojRozregNieznRozkl-Pi_n_Rekurencja})
we get
\small
{%\setlength\arraycolsep{1pt}
 \begin{eqnarray}
\label{PojRozregNieznRozkl-Pi_n-d-1_Od Pi_n-d-1}
\Pi_{n-d-1}^{i,j} &=&\left[(q_{ij} - \Pi_{n}^{i,j})\sum_{k=0}^d p_{ij}^{d-k}L_{k+1}^{i,j}(\underline{X}_{n-d-1,n}) - \Pi_{n}^{i,j}p_{ij}^{d+1}L_{0}^{i,j}(\underline{X}_{n-d-1,n})\right] \\
&& \times \Bigg[(1-\Pi^{i,j}_n)\Big( q_{ij}\sum_{k=0}^d p_{ij}^{d-k}L_{k+1}^{i,j}(\uX_{n-d-1,n})-L^{i,j}_{d+1}(\uX_{n-d-1,n})\Big) \nonumber\\
&&\;\; - \Pi^{i,j}_n p_{ij}^{d+1}L_{0}^{i,j}(\uX_{n-d-1,n})\Bigg]^{-1} \nonumber
\end{eqnarray}
}
\normalsize
The result (\ref{PojRozregNieznRozkl-Pi_n-d-1_Od Pi_n-d-1}) and formula (\ref{PojRozregNieznRozkl-DrugiCzlonFunkcjiWyplaty1}) lead us to:
\small
{%\setlength\arraycolsep{-1pt}
 \begin{eqnarray}
\label{PojRozregNieznRozkl-PodstWyplataJakoWWO-Drugi czlon}
\phiP&(& \theta \leq  n-d-1 \mid \cF_n, \ubeta = (i,j) ) \\
   &=&\;\; \frac{p_{ij}^{d+1}L_{0}^{i,j}(\underline{X}_{n-d-1,n}) \Pi_{n}^{i,j}- q_{ij}\sum_{k=0}^d p_{ij}^{d-k}L_{k+1}^{i,j}(\underline{X}_{n-d-1,n})(1- \Pi_{n}^{i,j})}{p_{ij}^{d+1}L_{0}^{i,j}(\underline{X}_{n-d-1,n}) }.\nonumber
\end{eqnarray}}
\normalsize
Applying equations (\ref{PojRozregNieznRozkl-PierwszyCzlonFunkcjiWyplaty}) and (\ref{PojRozregNieznRozkl-PodstWyplataJakoWWO-Drugi czlon})
in formula (\ref{PojRozregNieznRozkl-PodstWyplataJakoWWO1}) we get the thesis.
\end{pf}

Notice that for $n\geq d+1$ function $h$ under expectation in (\ref{PojRozregNieznRozkl-StaraWyplataNowaWyplata1}) depends
on process $\eta_n = (\underline{X}_{n-d-1,n},\ovPi_n, \ovBet_n)$. It turns out that $\{\eta_n\}$ is Markov Random Function
(see Lemma \ref{PojRozregNieznRozkl-FunkcjaMarkowska} in Appendix \ref{PojRozregNieznRozkl-dodatek}). We do not care about
$\{\eta_{n}\}$ for $n < d+1$. It is a consequence of discussion in Lemma \ref{PojRozregNieznRozkl-lematCzasStopuMax} which
leads to the conclusion that under the considered payoff function (criterion) it is not optimal to stop before instant
$d+1$. The decision maker can start his decision based on at least $d+1$ observations $X_1,\ldots, X_{d+1}$.

Lemmata \ref{PojRozregNieznRozkl-ZmianaWyplaty} and \ref{PojRozregNieznRozkl-FunkcjaMarkowska} imply that initial problem
can be reduced to the optimal stopping of Markov Random Function $(\eta_n, \cF_n, \bP_{\underline{y}}^\varphi)_{n=1}^\infty$,
where $\uy=(\ux_{n-d-1,n},\ovgamma,\ovdelta)\in{\itXi}= \bbE^{d+2}\times [0,1]^{l_1l_2}\times [0,1]^{l_1l_2}$ with payoff described by (\ref{PojRozregNieznRozkl-NowaWyplata}). However, the new problem is no longer homogeneous one as it is emphasized by the definition of $\uy$. It is a consequence of the fact that the process $\{\eta_n\}$ for $n< d+1$ has formally different structure than for $n\geq d+1$. Thus, the payoffs for instances $n \leq d+1$ are different. Lemma~\ref{PojRozregNieznRozkl-FunkcjaMarkowska} gives a justification to work with the homogeneous part of the process in construction the optimal estimator of the disorder moment.

To solve the maximization problem (\ref{PojRozregNieznRozkl-StaraWyplataNowaWyplata1}), for any Borel function
$u: {\itXi}\longrightarrow \R$ let us define operators:
\begin{eqnarray}
  \bT u(\underline{x}_{1,d+2},\overline{\gamma}, \overline{\delta}) &=& \phiE_{(\ux_{1,d+2},\ovgamma,\ovdelta)}\left[u(\underline{X}_{n-d,n+1},\ovPi_{n+1}, \ovBet_{n+1}) \right] \nonumber\\
   &=& \phiE\left[u(\underline{X}_{n-d,n+1},\ovPi_{n+1}, \ovBet_{n+1}) \mid  (\underline{X}_{n-d-1,n},\ovPi_{n}, \ovBet_{n}) = (\ux_{1,d+2},\ovgamma,\ovdelta)\right], \nonumber\\
  \bQ u(\underline{x}_{1,d+2},\overline{\gamma}, \overline{\delta}) &=& \max\{u(\underline{x}_{1,d+2},\overline{\gamma}, \overline{\delta}), \bT u(\underline{x}_{1,d+2},\overline{\gamma}, \overline{\delta}) \}.  \nonumber
\end{eqnarray}

Operators $\bT$ and $\bQ$ act on function $h$ and they determine the shape of optimal stopping rule $\tau^{\star}$.
Recursive formulas are given by Lemma \ref{PojRozregNieznRozkl-OperatoryDlaNowejWyplaty}, which is presented in
Appendix \ref{PojRozregNieznRozkl-dodatek}. Lemma \ref{PojRozregNieznRozkl-OperatoryDlaNowejWyplaty}
characterizes structure of sequence of functions $s_k(\underline{x}_{1,d+2},\overline{\gamma}, \overline{\delta})$, where $\underline{x} \in \bbE^{d+2}$, $\overline{\gamma}, \overline{\delta} \in [0,1]^{l_1l_2}$, which is used in the theorem stated below.
\begin{ptheorem}
    \label{PojRozregNieznRozkl-TwierdzPodsumowanie}
The solution of problem (\ref{PojRozregNieznRozkl-Problem}) is the following stopping rule:
{
\footnotesize
        \begin{eqnarray}
        \label{PojRozregNieznRozkl-optymalnyStop}
\!\!\!\!\!\!\!\!\!\!\!\!\!\!\!\!\!\!\! \tau^{*} &=&
            \left\{
    \begin{array}{ll}
                 \inf\left\{ \frac{}{}n\geq d+2: (\underline{X}_{n-1-d,n},\ovPi_{n}, \ovBet_{n}) \in D^{\star} \right\},& \mbox{if $\widetilde{h}(\ovPi_{d+1},\ovBet_{d+1}) < s^{*}({\uX}_{1,d+2},\ovPi_{d+2}, \ovBet_{d+2})$,}\\
                     &  \\
                d+1, &\mbox{if $\widetilde{h}(\ovPi_{d+1},\ovBet_{d+1}) \geq s^{*}({\uX}_{1,d+2},\ovPi_{d+2}, \ovBet_{d+2})$,}\\
            \end{array}
    \right.
\end{eqnarray}
\normalsize}
where the stopping area $D^{\star}$:
{
%\small
%\setlength\arraycolsep{0.3pt}
\begin{eqnarray*}
D^{\star} = \Bigg\{ &(&\underline{x}_{1,d+2},\overline{\gamma}, \overline{\delta}) \in {\itXi}: h(\underline{x}_{1,d+2},\overline{\gamma}, \overline{\delta})  \geq
s^{*}(\underline{x}_{1,d+2},\overline{\gamma}, \overline{\delta})
\Bigg\}, \nonumber
\end{eqnarray*}
\normalsize
}
and  $s^{*}({\ux}_{1,d+2},\overline{\gamma}, \overline{\delta}) = \lim_{k\longrightarrow \infty }s_k({\ux}_{1,d+2},\overline{\gamma}, \overline{\delta})$.
\end{ptheorem}
\begin{pf}
First let us consider subproblem of finding the optimal rule
$\widetilde{\tau}^{\star} \in \frF^X_{d+2}$:
\begin{eqnarray}
\label{PojRozregNieznRozkl-PrzeformulowanyProblem_dplus2}
       \phiE\left[ h(\underline{X}_{\widetilde{\tau}^{\star}-d-1,\widetilde{\tau}^{\star}},\ovPi_{\widetilde{\tau}^{\star}}, \ovBet_{\widetilde{\tau}^{\star}}) \right] = \sup_{\tau \in \frF^X_{d+2}}\phiE\left[ h(\underline{X}_{\tau-d-1,\tau},\ovPi_{\tau}, \ovBet_{\tau}) \right].
\end{eqnarray}
Then, basing on Lemmata \ref{PojRozregNieznRozkl-lematCzasStopu}, \ref{PojRozregNieznRozkl-lematCzasStopuMax} and according
to optimal stopping theory (see \cite{shi78:optimal}) it is known that $\tau_0$ defined by
(\ref{PojRozregNieznRozkl-StopIntuicyjny}) can be expressed as
\[
    \tau_0 = \inf\{n \geq d+2: h(\underline{X}_{n-1-d,n},\ovPi_{n}, \ovBet_{n}) \geq h^{*}(\underline{X}_{n-1-d,n},\ovPi_{n}, \ovBet_{n}) \},
\]
where $h^{*}(\underline{x}_{1,d+2},\overline{\gamma}, \overline{\delta}) = \lim_{k \longrightarrow \infty} \bQ^{k}h(\underline{x}_{1,d+2},\overline{\gamma}, \overline{\delta})$. The limit exists according to the Lebesgue's theorem and structure of functions $h$ and $s_k$. Lemma \ref{PojRozregNieznRozkl-OperatoryDlaNowejWyplaty} implies that:
{
%\footnotesize
%\setlength\arraycolsep{1pt}
\begin{eqnarray*}
\tau_0 &=& \inf\left\{n \geq d+2: h({\uX}_{n-1-d,n},\ovPi_{n}, \ovBet_{n})\right.\\
\nonumber&&\mbox{}\geq \left.\max\left\{h({\uX}_{n-1-d,n},\ovPi_{n}, \ovBet_{n}),  s^{*}({\uX}_{n-1-d,n},\ovPi_{n}, \ovBet_{n})\right\} \right\} \nonumber\\
&=& \inf \left\{n \geq d+2: h(\underline{X}_{n-1-d,n},\ovPi_{n}, \ovBet_{n})  \geq s^{*}(\underline{X}_{n-1-d,n},\ovPi_{n}, \ovBet_{n})\right\}. \nonumber
\end{eqnarray*}
%\normalsize
}
%The next step is to obtain the value of subproblem defined for $n \geq d+2$. We have:
%{
%\small
%\setlength\arraycolsep{0pt}
%\begin{eqnarray*}
%\phiE\big( h^{*}(\underline{X}_{1,d+2},\overline{\Pi}_{d+2}, \overline{B}_{d+2})\big)&=& \phiE \left( \max\left\{h(X_{1,d+2},\Pi^{i,j}_{d+2})B^{i,j}_{d+2}),  s^{*}(\underline{X}_{1,d+2},\ovPi_{d+2}, \ovBet_{d+2})\right\} \right) \nonumber \\
%&\hspace{-24em}=&\hspace{-12em} \int_{\bbE^{d+2}}\max\left\{h({\ux}_{1,d+2},%
%\itPi^{i,j}_{d+2},\itBet_{d+2}),   s^{*}(\underline{x}_{1,d+2},\itPi_{d+2}, \itBet_{d+2})\right\} \nonumber \\
%&& \times \big[ \sum_{i,j}\itPsi^{i,j}({d+2},\ux_{0,d+2},\pi_{ij})b_{ij} \big] \mu(d\ux_{0,d+2})=  V \nonumber
%\end{eqnarray*}
%\normalsize}
%where $\itPi_{d+2}$, $\itBet_{d+2}$ are functions of $\ux_{1,d+2},\ovp,\ovbet$.
%%%%% TUTAJ Zacznij %%%%%%%%%%%%%%%%%%%%
%The third equality results from (\ref{PojRozregNieznRozkl-Transformacja_na_potrzeby_rekurencji Pin8_2}) where we put $n=d+2$ and $l=d+1$.

According to optimality principle rule $\widetilde{\tau}^{*}$ solves maximization problem of
(\ref{PojRozregNieznRozkl-StaraWyplataNowaWyplata1}) if only at $n = d+1$ the payoff $\widetilde{h}$ will be smaller than
expected payoff in successive periods (for $n > d+1$). Thus, another words:
{
\begin{eqnarray}
\tau^{\star} = \widetilde{\tau}^{*},\mbox{ \quad if $\widetilde{h}(\uX_{0,d+1},\ovPi_{d+1},\ovBet_{d+1}) < s^{*}({\uX}_{1,d+2},\ovPi_{d+2}, \ovBet_{d+2})$.}
\end{eqnarray}}
In opposite case $\tau^{\star} = d+1$. This ends the proof of formula (\ref{PojRozregNieznRozkl-optymalnyStop}).
\end{pf}
%\begin{premark}
%Computing $\tau^{\star}$ is connected with calculation of value $V$ at instant $n=d+1$. Thus we have to know
%$\ovPi_{d+2}$ i $\ovBet_{d+2}$. Those vectors can be computed for each realization of
%$\underline{X}_{1,d+2} = \underline{x}_{1,d+2}$ basing on recursive formulas (\ref{PojRozregNieznRozkl-Pi_n_Rekurencja2})
%and (\ref{PojRozregNieznRozkl-B_n_Rekurencja}) obtained in lemmas \ref{PojRozregNieznRozkl-dodatekWzor4} and
%\ref{PojRozregNieznRozkl-dodatekWzor5}
%\end{premark}

%
\subsection{Acknowledgements}
We have benefited from the remarks of anonymous referee of submitted presentation to the Program Committee of IWSM 2009. He has provided numerous corrections to the manuscript.

%\input{UnDistSinglAppendix_1.tex}
%%%%%%%%%%%%%%%%%%%%%%%%%%%%%%%%%%% Unspecified distribution in single disorder problem-Appendix
\appendix
\section{Lemmata}\label{PojRozregNieznRozkl-dodatek}
In appendix we present useful formulae and lemmata which help to obtain solution of problem (\ref{PojRozregNieznRozkl-Problem}).
\begin{premark}\label{PojRozregNieznRozkl-dodatekWzor1}
For $n \geq l \geq 0$, $k >0$, $i \in I_0$, $j \in I_1$, on the set \samepage $\{\omega: \underline{X}_{0,l} \in \underline{A}_{0,l}, A_0 = \{x\}, A_i \in \cF_i, i \leq l\}$
the following equations hold:
{
\setlength\arraycolsep{0pt}
\begin{eqnarray}
    \label{PojRozregNieznRozkl-pq1}
\phiP&(& \theta = n+k \mid \underline{X}_{0,l} \in \underline{A}_{0,l}, \ubeta = (i,j), \theta > n) \nonumber\\
&=&\;\;
\left\{
    \begin{array}{ll}
                p_{ij}^{k-1}q_{ij},  & \mbox{if $n, k > 0$}, \\
                \pi_{ij},& \mbox{if $n = 0, k = 1$,}\\
                (1-\pi_{ij})p_{ij}^{k-2}q_{ij}, \;\;\; & \mbox{if $n = 0, k > 1$,} \\
            \end{array}
    \right. \nonumber
\end{eqnarray}}
\end{premark}

\begin{premark}
\begin{enumerate}
\item The simple consequence of the formula (\ref{PojRozregNieznRozkl-pq1}) we get
\begin{eqnarray}
    \label{PojRozregNieznRozkl-pq2}
\phiP&(& \theta > n+k \mid \underline{X}_{0,l} \in \underline{A}_{0,l}, \ubeta = (i,j), \theta > n) \nonumber\\
&=&\;\;\left\{
    \begin{array}{ll}
                p_{ij}^{k},  & \mbox{if $n, k > 0$}, \\
                (1-\pi_{ij})p_{ij}^{k-1},\;\;\;  & \mbox{if $n = 0, k > 0$.} \\
            \end{array}
    \right.
\end{eqnarray}
\item Formula (\ref{PojRozregNieznRozkl-pq2}) for $k=1$ is given by:
{
\setlength\arraycolsep{0pt}
\begin{eqnarray}
   \label{PojRozregNieznRozkl-pq3}
\phiP&(& \theta \neq n+1 \mid \underline{X}_{0,l} \in \underline{A}_{0,l}, \ubeta = (i,j), \theta > n) \nonumber\\
 &=& \phiP( \theta > n+1 \mid \underline{X}_{0,l} \in \underline{A}_{0,l}, \ubeta = (i,j), \theta > n) \nonumber\\
 &=& \left\{
    \begin{array}{ll}
                p_{ij},  & \mbox{if $n >0$}, \\
                (1-\pi_{ij}),\;\;\;  & \mbox{if $n = 0$.} \\
            \end{array}
    \right.
\end{eqnarray}
}
\end{enumerate}
\end{premark}
%%%Equations (\ref{PojRozregNieznRozkl-pq1})-(\ref{PojRozregNieznRozkl-pq3}) will be used in proofs of another lemmas
%%%collected in appendix.

\begin{lemma}\label{PojRozregNieznRozkl-dodatekWzor2}
For $n > 0$, $k \geq 0$, $i \in I_0$, $j \in I_1$ the following equation is satisfied:
\begin{eqnarray}
    \label{PojRozregNieznRozkl-PierwszyCzlonFunkcjiWyplaty}
\phiP( \theta \leq n+k \mid \tilde{\cF}_n ) =  1 -p_{ij}^k(1-\Pi_n^{i,j}).
\end{eqnarray}
\end{lemma}
\begin{pf} We are going to show equality on the set
$\{\omega: \underline{X}_{0,n} \in \underline{B}_{0,n}, B_0 = \{x\}\}$
{
\setlength\arraycolsep{0pt}
\begin{eqnarray}\label{phiP1}
\phiP( \theta > n+k, \underline{X}_{0,n} \in \underline{B}_{0,n}, \ubeta = (i,j) ) &=&\hspace{-2em} \displaystyle{\int\limits_{\{\omega:\uX_{0,n}\in\underline{B}_{0,n},\ubeta=(i,j)\}}}\hspace{-2em}\one_{\{\omega:\theta>n+k\}}d\phiP \\
\nonumber&\hspace{-24em}=&\hspace{-12em}\hspace{-2em} \displaystyle{\int\limits_{\{\omega:\uX_{0,n}\in\underline{B}_{0,n},\ubeta=(i,j)\}}}\hspace{-2em}\phiE(\one_{\{\omega:\theta>n+k\}}|\tilde{\cF}_n)d\phiP\\
\nonumber&\hspace{-24em}=&\hspace{-12em}\hspace{-2em}\displaystyle{\int\limits_{\{\omega:\uX_{0,n}\in\underline{B}_{0,n}\}}}\hspace{-2em}\phiE(\one_{\{\omega:\ubeta=(i,j)\}}\phiE(\one_{\{\omega:\theta>n+k\}}|\tilde{\cF}_n)|\cF_n)d\phiP.
\end{eqnarray}
}
By direct computation calculation we get
{
\setlength\arraycolsep{0pt}
\begin{eqnarray}\label{phiP2}
\phiP( \theta > n+k, \underline{X}_{0,n} \in \underline{B}_{0,n}, \ubeta = (i,j) )\\
\nonumber&\hspace{-24em}=&\hspace{-12em} \int_{\underline{B}_{0,n}}\sum_{s=n+k+1}^\infty (1-\pi_{ij})q_{ij}b_{ij}p_{ij}^{s-2}L_0^{ij}(\ux_{0,n})\mu(d\ux_{0,n})\\
\nonumber&\hspace{-24em}=&\hspace{-12em}  p_{ij}^{k}\int_{\underline{B}_{0,n}}\frac{(1-\pi_{ij})b_{ij}p_{ij}^{n-1}L_0^{ij}(\ux_{0,n})}{S_n^{i,j}(\ux_{0,n})}\frac{S_n^{i,j}(\ux_{0,n})}{S_n(\ux_{0,n})}S_n(\ux_{0,n})\mu(d\ux_{0,n})\\
\nonumber&\hspace{-24em}=&\hspace{-12em}p_{ij}^{k}\hspace{-2em} \displaystyle{\int\limits_{\{\omega:\uX_{0,n}\in\underline{B}_{0,n}\}}}\hspace{-2em}\one_{\{\omega:\theta>n\}}\one_{\{\omega:\ubeta=(i,j)\}}d\phiP =p_{ij}^{k}\hspace{-2em} \displaystyle{\int\limits_{\{\omega:\uX_{0,n}\in\underline{B}_{0,n}\}}}\hspace{-2em}(1-\Pi^{i,j}_n)B_n^{i,j}d\phiP\\
&\hspace{-24em}=&\hspace{-12em}p_{ij}^{k}\hspace{-2em}\displaystyle{\int\limits_{\{\omega:\uX_{0,n}\in\underline{B}_{0,n}\}}}\hspace{-2em}\phiP(\ubeta=(i,j)|\cF_n)\phiP(\theta> n|\tilde{\cF}_n)d\phiP
\end{eqnarray}
}
Henceforth we have
\[
\phiE(\one_{\{\omega:\ubeta=(i,j)\}}\phiE(\one_{\{\omega:\theta>n+k\}}|\tilde{\cF}_n)|\cF_n)=p_{ij}^{k}\phiP(\ubeta=(i,j)|\cF_n)\phiP(\theta> n|\tilde{\cF}_n).
\]
Comparison of (\ref{phiP1}) and (\ref{phiP2}) implies \ref{PojRozregNieznRozkl-PierwszyCzlonFunkcjiWyplaty} and this ends the proof of lemma.
%In above chain of reasoning the result of lemma \ref{PojRozregNieznRozkl-dodatekWzor1} has been used.

\end{pf}
\begin{lemma}\label{PojRozregNieznRozkl-dodatekWzor3}
For $n >k \geq 0$, $i \in I_0$, $j \in I_1$ it is true that
\small
\begin{eqnarray}
    \label{PojRozregNieznRozkl-DrugiCzlonFunkcjiWyplaty1}
    \phiP( \theta \leq n-k-1 \mid \tilde{\cF}_n) =1-(1-\Pi_n^{i,j})\left(1+q_{ij}\sum_{s=1}^{k+1}\frac{L_s^{i,j}(\underline{X}_{n-s+1,n})}{p_{ij}^sL_0^{i,j}(\underline{X}_{n-s+1,n}}\right).
\end{eqnarray}
\normalsize
\end{lemma}
\begin{pf} If $n=k+1$ then
\begin{eqnarray}
\phiP( \theta \leq n-k-1 \mid \tilde{\cF}_n) = \phiP( \theta \leq 0 \mid \tilde{\cF}_{k+1}) = 0. \nonumber
\end{eqnarray}
Because of the fact that $\theta >0$ a.s.:
\begin{eqnarray}
\label{Pi0}
\Pi_{n-k-1}^{i,j} = \Pi_0^{i,j} = \phiP( \theta \leq 0 \mid \tilde{\cF}_{0}) = 0.
\end{eqnarray}
Hence formula (\ref{PojRozregNieznRozkl-DrugiCzlonFunkcjiWyplaty1}) holds. The case where $n > k+1$
we have
\begin{equation*}
\phiP( \theta > n-k-1 \mid \tilde{\cF}_n)=\phiP( \theta >n \mid \tilde{\cF}_n)+\sum_{s=1}^{k+1}\phiP( \theta= n-s \mid \tilde{\cF}_n).
\end{equation*}
On the set $D_n=\{\omega: \ubeta=(i,j), \underline{X}_{0,n} \in \underline{B}_{0,n}, B_0 = \{x\}\}$ we have
\begin{eqnarray*}
\phiP( \theta= n-s, D_n)&=&\int_{D_n}\one_{\{\theta=n-s\}}d\phiP =\int_{D_n}\phiP(\theta=n-s|\tilde{\cF}_n) d\phiP \\
&\hspace{-10em}=&\hspace{-5em}\hspace{-1em} \displaystyle{\int\limits_{\times_{r=1}^n B_r}}\hspace{-1em}(1-\pi_{ij})p_{ij}^{n-s-2}q_{ij}b_{ij}L_{s+1}^{i,j}(\underline{x}_{n-s,n})d\mu(\underline{x}_{0,n}) \\
&\hspace{-10em}=&\hspace{-5em}\hspace{-1em} \displaystyle{\int\limits_{\times_{r=1}^n B_r}}\hspace{-1em}\frac{q_{ij}L_{s+1}^{i,j}(\underline{x}_{n-s,n})}{p_{ij}^{s+1}L_{0}^{i,j}(\underline{x}_{n-s,n})}\frac{(1-\pi_{ij})p_{ij}^{n-1}b_{ij}L_{0}^{i,j}(\underline{x}_{0,n})}{S_n^{i,j}(\underline{x}_{0,n})}\frac{S_n^{i,j}(\underline{x}_{0,n})}{S_n(\underline{x}_{0,n})}S_n(\underline{x}_{0,n})d\mu(\underline{x}_{0,n})\\
&\hspace{-10em}=&\hspace{-5em}p_{ij}^{k}\hspace{-1  em} \displaystyle{\int\limits_{\{\omega:\uX_{0,n}\in\underline{B}_{0,n}\}}}\hspace{-1em}\frac{q_{ij}L_{s+1}^{i,j}(\underline{X}_{n-s,n})}{p_{ij}^{s+1}L_{0}^{i,j}(\underline{X}_{n-s,n})}(1-\Pi_n^{ij})B_n^{i,j}d\phiP.
\end{eqnarray*}
Therefore
\begin{equation*}
\phiP(\theta=n-s|\tilde{\cF}_n)=\frac{q_{ij}L_{s+1}^{i,j}(\underline{X}_{n-s,n})}{p_{ij}^{s+1}L_{0}^{i,j}(\underline{X}_{n-s,n})}(1-\Pi^{i,j}_n)
\end{equation*}
and
\begin{equation*}
\phiP(\theta>n-k-1|\tilde{\cF}_n)=\left(1+q_{ij}\sum_{s=0}^{k}\frac{L_{s+1}^{i,j}(\underline{X}_{n-s,n})}{p_{ij}^{s+1}L_{0}^{i,j}(\underline{X}_{n-s,n})}\right)(1-\Pi^{i,j}_n).
\end{equation*}
\end{pf}

\begin{lemma}\label{PojRozregNieznRozkl-dodatekWzor4}
For $n >l \geq 0$, $i \in I_0$, $j \in I_1$ following equation holds:
{
\setlength\arraycolsep{0pt}
 \begin{eqnarray}
\label{PojRozregNieznRozkl-Pi_n_Rekurencja}
\Pi_{n}^{i,j}&=& \left\{
    \begin{array}{ll}
\Pi^{i,j}(l,\uX_{n-l-1,n},\Pi_{n-l-1}^{i,j}),& \mbox{if $n > l+1$}, \\ \\
                     %& \mbox{if $n > l+1$}, \\
\frac{\itLambdaT(l,\uX_{n-l-1,n},\pi_{ij})}{\itPsiT^{i,j}(l,\uX_{0,l+1},\pi_{ij})}, & \mbox{if $n = l+1$.} \\
            \end{array}
    \right.
\end{eqnarray}
\normalsize
}
\end{lemma}

\begin{premark}
In particular, taking $l=0$, we get equation characterizing "one-step" dynamics of the process $\Pi_n^{i,j}$:
 \begin{eqnarray}
\label{PojRozregNieznRozkl-Pi_n_Rekurencja2}
\Pi_n^{i,j}
&=& \left\{
    \begin{array}{ll}
                \frac{ f_{X_{n-1}}^{1,j}(X_n)(q_{ij}+p_{ij}\Pi_{n-1}^{i,j}) }{f_{X_{n-1}}^{1,j}(X_n)(q_{ij}+p_{ij}\Pi_{n-1}^{i,j}) + f_{X_{n-1}}^{0,i}(X_n)p_{ij}(1-\Pi_{n-1}^{i,j})},& \mbox{if $n > 1$}, \\
                     %& \mbox{if $n > l+1$}, \\
                \frac{f_{X_{0}}^{1,j}(X_1)\pi_{ij}}{f_{X_{0}}^{1,j}(X_1)\pi_{ij}+f_{X_{0}}^{0,i}(X_1)(1-\pi_{ij})}, & \mbox{if $n = 1$,} \\
            \end{array}
    \right.
\end{eqnarray}
with initial condition $\Pi_0^{i,j} = 0$.
\end{premark}
\begin{premark}
For $l>0$ recursive structure defined in equation (\ref{PojRozregNieznRozkl-Pi_n_Rekurencja}) requires vector of initial
states $\Pi_0^{i,j}, \Pi_1^{i,j}, \ldots, \Pi_l^{i,j}$. State $\Pi_0^{i,j}$ is given above. To obtain remaining states
$\Pi_1^{i,j}, \ldots, \Pi_l^{i,j}$ it is enough to apply formula (\ref{PojRozregNieznRozkl-Pi_n_Rekurencja2}).
\end{premark}
%\end{lemma}
\begin{pf} (of Lemma~\ref{PojRozregNieznRozkl-dodatekWzor4}) Condition $\Pi_0^{ij}=0$ has been shown in lemma \ref{PojRozregNieznRozkl-dodatekWzor3} (equation (\ref{Pi0})).
We have following recursive relation;
\small
\setlength\arraycolsep{0pt}
\begin{eqnarray*}
S_{n-s-1}^{i,j}(\ux_{0,n-s-1})\itPsi^{i,j}(\ux_{n-s-1,n},\Pi^{i,j}(n-s,\ux_{n-s-1,n},\pi_{ij}))&\mbox{}& \nonumber\\
&\hspace{-45em}=&\hspace{-22em}S_{n-s-1}^{i,j}(\ux_{0,n-s-1}) \Pi_{n-s-1}^{i,j} L_{l+1}(\ux_{n-s-1,n})+ S_{n-s-1}(\ux_{0,n-s-1})(1-\Pi_{n-s-1}^{i,j})  \nonumber\\
&\hspace{-40em}&\hspace{-19em} \times \left[q_{ij} \sum_{k=0}^{s}p_{ij}^{s-k} L_{k+1}(\ux_{n-s-1,n}) + p_{ij}^{s+1}L_0(\ux_{n-s-1,n})\right] \nonumber\\
&\hspace{-45em}=&\hspace{-22em}\left(\pi_{ij}L_{n-s}^{i,j}(\ux_{0,n-s-1})+(1-\pi_{ij})q_{ij}\sum_{k=1}^{n-s-1}p_{ij}^{k-1} L_{n-s-k}(\ux_{0,n-s-1})\right)L_{s+1}(\ux_{n-s-1,n})\\
&\hspace{-40em}&\hspace{-19em}  + (1-\pi_{ij})p_{ij}^{n-s-1}L_0(\ux_{0,n-s-1})\left( \sum_{k=0}^{l}p_{ij}^{s-k}q_{ij} L_{k+1}(\ux_{n-s-1,n}) + p_{ij}^{s+1}L_0(\ux_{n-s-1,n})\right)\nonumber\\
\nonumber&\hspace{-45em}=&\hspace{-22em}\pi_{ij}L_{n}^{i,j}(\ux_{0,n})+(1-\pi_{ij}\Big[\sum_{k=1}^{n-s-1}p_{ij}^{k-1}q_{ij} L_{n-k+1}(\uX_{0,n})\\
&\hspace{-40em}&\hspace{-19em}+ \sum_{k=0}^{s}p_{ij}^{n-k-1}q_{ij} L_{k+1}(\ux_{0,n}) + p_{ij}^{n}L_0(\uX_{0,n})\Big]\nonumber\\
&\hspace{-45em}=&\hspace{-22em}\pi_{ij}L_{n}^{i,j}(\ux_{0,n})+(1-\pi_{ij})\Big[\sum_{k=1}^{n-s-1}p_{ij}^{k-1}q_{ij} L_{n-k+1}(\uX_{0,n})\nonumber\\
&\hspace{-40em}&\hspace{-19em} + \sum_{k=n-s}^{n}p_{ij}^{k-1}q_{ij} L_{n-k+1}(\uX_{0,n}) + p_{ij}^{n}L_0(\uX_{0,n})\Big]\nonumber\\
&\hspace{-45em}=&\hspace{-22em}\pi_{ij}L_{n}^{i,j}(\ux_{0,n})+(1-\pi_{ij})\left[\sum_{k=1}^{n}p_{ij}^{k-1}q_{ij} L_{n-k+1}(\uX_{0,n})+ p_{ij}^{n}L_0(\uX_{0,n})\right]\nonumber\\
&\hspace{-45em}=&\hspace{-22em} S_n^{i,j}(\ux_{0,n}).\nonumber
\end{eqnarray*}
Now, on the set
$D_n=\{\omega:\underline{X}_{0,n}\in \underline{B}_{0,n}\}$, $X_0=x$ and $B_i\in\cB$ we have by (\ref{SnijFormula1}):
\begin{eqnarray*}
\bP(\theta>n,\ubeta=(i,j),D_n)&=&\hspace{-1em}\int\limits_{\{\ubeta=(i,j),D_n\}}\hspace{-1em}\one_{\{\theta>n\}}d\phiP
=\hspace{-1em}\int\limits_{\{\ubeta=(i,j),D_n\}}\hspace{-1em}\phiP(\theta>n|\tilde{\cF}_n)d\phiP\\
&\hspace{-20em}=&\hspace{-10em}  \int\limits_{\underline{B}_{0,n}}\frac{p_{ij}^{s-1}L_0^{ij}(\ux_{n-s-1,n})}{\itPsi^{i,j}(n-s,\ux_{n-s-1,n},\Pi^{i,j}(n-s,\ux_{n-s-1,n},\pi_{ij}))}\\
&\hspace{-20em}&\hspace{-10em}\times\frac{(1-\pi_{ij})p_{ij}^{n-s-2}L_0^{i,j}(\ux_{0,n-s-1})}{S_{n-s-1}^{i,j}(\ux_{0,n-s-1})}\frac{b_{ij}S_{n-s-1}^{i,j}(\ux_{0,n-s-1})}{S_n(\ux_{0,n})}S_n(\ux_{0,n})\mu(d\ux_{1,n})\\
&\hspace{-20em}=&\hspace{-10em}\int\limits_{D_n}\frac{p_{ij}^{s-1}L_0^{ij}(\uX_{n-s-1,n})}{\itPsi^{i,j}_{n-s-1}(\uX_{n-s-1,n},\Pi_{n-s-1}^{i,j})}(1-\Pi_{n-s-1}^{i,j})B_n^{i,j}d\phiP.
\end{eqnarray*}
This follows
\begin{equation*}
\phiP(\theta>n|\tilde{\cF}_n)=\frac{p_{ij}^{s-1}L_0^{ij}(\uX_{n-s-1,n})}{\Psi^{i,j}_{n-s-1}(\uX_{n-s-1,n},\Pi_{n-s-1}^{i,j})}(1-\Pi_{n-s-1}^{i,j}).
\end{equation*}
In the case where $n =s+1$ the proof is similar.
\end{pf}

\begin{lemma}\label{PojRozregNieznRozkl-dodatekWzor5}
For $n > 0  $, $i \in I_0$, $j \in I_1$ we have
{\setlength\arraycolsep{0pt}
 \begin{eqnarray}
\label{PojRozregNieznRozkl-B_n_Rekurencja}
B_n^{i,j} &=&\left\{
    \begin{array}{ll}
       \Gamma^{i,j}(0,\uX_{n-1,n},\ovBet_{n-1}, \ovPi_{n-1}),\;\; & \mbox{if $n >1$},\\
        &\\
        \frac{ b_{n-1}^{i,j}{\itPsiT}^{i,j}(0,\uX_{0,1},\pi_{ij}) }{\ST(0,\uX_{0,1},\ovb,\ovpi)}, &  \mbox{if $n = 1$}.\\
    \end{array}
            \right.
\end{eqnarray}}
with condition $B_0^{i,j} = b_{ij}$.
\end{lemma}
\begin{pf} First, let us verify the initial condition:
 \begin{eqnarray}
\label{PojRozregNieznRozkl-B_n_WarPocz}
B_0^{i,j} = \phiP(\ubeta = (i,j) \mid \cF_0) = \phiP(\ubeta = (i,j)) = b_{ij}\nonumber.
\end{eqnarray}
Let $n > 1$. Let us consider formula (\ref{PojRozregNieznRozkl-B_n_Rekurencja}) on the set
$D_n=\{\omega: \underline{X}_{0,n} \in \underline{B}_{0,n}; B_0 = \{x\}$,  $B_i\in\cB$ for $1\leq i \leq n\}$:
{\setlength\arraycolsep{0pt}
\begin{eqnarray}
\label{PojRozregNieznRozkl-B_n_Rekurencja1}
\!\!\!\!\!\!\!\!\!\!\phiP( \ubeta = (i,j),\underline{X}_{0,n} \in \underline{B}_{0,n} )&=&\int_{D_n}\one_{\{\ubeta=(ij)\}} d\phiP=\int_{D_n}\phiE(\one_{\{\ubeta=(ij)\}}|\cF_n) d\phiP \nonumber\\
 \!\!\!\!\!\!\!\!\!\!&=&\int_{\underline{B}_{1,n}} \frac{b_{ij}S_n^{i,j}(\ux_{0,n})}{S_n(\ux_{0,n})}S_n(\ux_{0,n})\mu(d\ux_{1,n})\nonumber \\
 \!\!\!\!\!\!\!\!\!\!&=&\int_{D_n} \frac{b_{ij}S_n^{i,j}(\uX_{0,n})}{S_n(\uX_{0,n})}d\phiP.
\end{eqnarray}}
Taking into account the formulae (\ref{Snij}), (\ref{denomAuxFunctPi}), (\ref{multidist}) and (\ref{auxFunctGamma}) we have gotten (\ref{PojRozregNieznRozkl-B_n_Rekurencja}) for $n>1$. The case $n=1$ is a consequence of (\ref{PojRozregNieznRozkl-B_n_Rekurencja1}) and (\ref{auxFunctGamma}) with (\ref{denomAuxFunctPiTilda}) and (\ref{multidistTilda}).
\end{pf}

\begin{lemma}
\label{PojRozregNieznRozkl-FunkcjaMarkowska}
Let $\eta_n = (\underline{X}_{n-d-1,n},\ovPi_n, \ovBet_n)$, where $n \geq d+1$. System $(\eta_n, \cF_n, \bP_{\underline{y}}^\varphi)$ is Markov Random Function.
\end{lemma}
\begin{pf} It is enough to show that $\eta_{n+1}$ is a function of $\eta_{n}$ and variable $X_{n+1}$ as well as that
conditional distribution of $X_{n+1}$ given $\cF_n$ depends only on $\eta_n$ (see \cite{shi78:optimal}).

For $x_1,...,x_{d+2}, y \in \bbE, \gamma_{ij}, \delta_{ij} \in [0,1]$, $i \in I_0, j \in I_1$
let us consider the following function
{
\small
 \begin{eqnarray}
\label{PojRozregNieznRozkl-FunkcjaMarkowskaRekur}
    \varphi&(&\underline{x}_{1,d+2},\overline{\gamma}, \overline{\delta},y)\nonumber \\
   &=& \left(\underline{x}_{2,d+2},y,\Pi^{1,1}(0,x_{d+2},y,\delta_{11}), \ldots, \Pi^{1,l_2}(0,x_{d+2},y,\delta_{1l_2}),\ldots,\right.  \nonumber \\
    && \;\;\;\Pi^{l_1,1}(0,x_{d+2},y,\delta_{l_11}), \ldots, \Pi^{l_1,l_2}(0,x_{d+2},y,\delta_{l_1l_2}), \nonumber \\
    && \;\;\;\Gamma^{1,1}(0,x_{d+2},y,\ovgamma, \ovdelta), \ldots, \Gamma^{1,l_2}(0,x_{d+2},y,\ovgamma, \ovdelta),\ldots, \nonumber\\
    && \;\;\;\left. \Gamma^{l_1,1}(0,x_{d+2},y,\ovgamma, \ovdelta), \ldots, \Gamma^{l_1,l_2}(0,x_{d+2},y,\ovgamma, \ovdelta) \right). \nonumber
 \end{eqnarray}
\normalsize
}
We will show that $\eta_{n+1} = \varphi(\eta_n,X_{n+1})$. Using formulas (\ref{PojRozregNieznRozkl-Pi_n_Rekurencja2})
and (\ref{PojRozregNieznRozkl-B_n_Rekurencja}) we express $\Pi_{n+1}^{i,j}$ as a function of $\Pi_n^{i,j}$
and $B_{n+1}^{i,j}$ as a function $B_n^{i,j}$. Then:
{
\small
\setlength\arraycolsep{0.3pt}
 \begin{eqnarray}
\label{PojRozregNieznRozkl-FunkcjaMarkowskaRekur1}
\varphi&(&\eta_{n}, X_{n+1}) \nonumber \\
    &=&\varphi(\underline{X}_{n-d-1,n},\ovPi_n, \ovBet_n, X_{n+1}) \nonumber \\
   &=& \left(\underline{X}_{n-d,n},X_{n+1},\Pi^{1,1}(0,\uX_{n,n+1},\Pi_n^{1,1}), \ldots, \Pi^{1,l_2}(0,\uX_{n,n+1},\Pi_n^{1,l_2}),\ldots,\right.  \nonumber \\
    && \;\;\;\Pi^{l_1,1}(0,\uX_{n,n+1},\Pi_n^{l_1,1}), \ldots, \Pi^{l_1,l_2}(0,\uX_{n,n+1},\Pi_n^{l_1,l_2}), \nonumber \\
    && \;\;\;\Gamma^{1,1}(0,\uX_{n,n+1},\ovBet_{n}, \ovPi_{n}), \ldots, \Gamma^{1,l_2}(0,\uX_{n,n+1},\ovBet_{n}, \ovPi_{n}),\ldots, \nonumber\\
    && \;\;\;\left. \Gamma^{l_1,1}(0,\uX_{n,n+1},\ovBet_{n}, \ovPi_{n}), \ldots, \Gamma^{l_1,l_2}(0,\uX_{n,n+1},\ovBet_{n}, \ovPi_{n}) \right). \nonumber\\
   &=& (\underline{X}_{n-d,n+1},\ovPi_{n+1}, \ovBet_{n+1}) = \eta_{n+1}. \nonumber
 \end{eqnarray}
\normalsize
}
Let us consider now the conditional expectation $u(X_{n+1})$ under the condition of $\sigma$~-field $\cF_n$, for Borel function $u: \bbE\longrightarrow \Re$. Applying equation (\ref{PojRozregNieznRozkl-PierwszyCzlonFunkcjiWyplaty}) ($k=1$) we get:
{
%\small
%\setlength\arraycolsep{0.3pt}
\begin{eqnarray}\label{expVALUE}
\phiE(u(X_{n+1})\mid \cF_n)&=&\sum_{i,j}\phiE(u(X_{n+1})\one_{\{\ubeta=(i,j)\}}\mid \cF_n) \\
\nonumber&=&\sum_{i,j}\phiE(u(X_{n+1})\one_{\{\theta\leq n+1\}}\one_{\{\ubeta=(i,j)\}}\mid \cF_n)\\
&&+\sum_{i,j}\phiE(u(X_{n+1})\one_{\{\theta> n+1\}}\one_{\{\ubeta=(i,j)\}}\mid \cF_n) \nonumber\\
\nonumber &=&\sum_{i,j} B_n^{i,j}\left[\phiE\left(p_{ij}\int_\bbE u(y)(1-\Pi^{i,j}(0,y,\Pi_n^{i,j}))f^{0,i}_{X_n}(y)\mu(dy) \mid \cF_n\right)\right.\\
       &&\left. + \phiE\left(\int_\bbE u(y)(q_{ij}+p_{ij}\Pi^{i,j}(0,y,\Pi_n^{i,j}))f^{1,j}_{X_n}(X_{n+1})\mu(dy) \mid \cF_n\right)\right] \nonumber
\end{eqnarray}
\normalsize
}
We see that conditional distribution of $X_{n+1}$ given $\cF_n$ depends only on component of $\eta_n$ what ends the proof.

\end{pf}

\begin{lemma}
 \label{PojRozregNieznRozkl-OperatoryDlaNowejWyplaty}
    Let
    {
\small
\setlength\arraycolsep{1pt}
\begin{eqnarray}
\label{PojRozregNieznRozkl-def_rk}
  s_k(\underline{x}_{1,d+2},\overline{\gamma}, \overline{\delta}) &=&
        \left\{
            \begin{array}{ll}
               \bT\;\bQ^k h(\underline{x}_{1,d+2},\overline{\gamma}, \overline{\delta}), & \mbox{if $k \geq 1$}, \\
               \bT\;h(\underline{x}_{1,d+2},\overline{\gamma}), & \mbox{if $k=0$.}\\
            \end{array}
    \right.
  %  s_k(\underline{x}_{1,d+2},\overline{\gamma}, \overline{\delta}) &=& \bT\;\bQ^k h(\underline{x}_{1,d+2},\overline{\gamma}, \overline{\delta}), \;k \geq 1, \nonumber\\
  %  s_0(\underline{x}_{1,d+2},\overline{\gamma}, \overline{\delta}) &=& \bT\;h(\underline{x}_{1,d+2},\overline{\gamma}), \;k = 0. \nonumber
\end{eqnarray}
\normalsize
}
    Then, for function $h(\underline{x}_{1,d+2},\overline{\gamma}, \overline{\delta})$ given by
    (\ref{PojRozregNieznRozkl-NowaWyplata}) and $k \geq 1$, following equalities hold:
{
\small
\setlength\arraycolsep{0pt}
\begin{eqnarray}
 \bQ^{k} h&(&\underline{x}_{1,d+2},\overline{\gamma}, \overline{\delta}) \nonumber\\
  &=& \;\max\left\{\sum_{i,j}\left( 1-p_{ij}^d+ q_{ij}\sum_{m=1}^{d+1}\frac{L_{m}^{i,j}(\underline{x}_{1,d+2})} { p_{ij}^m L_{0}^{i,j}(\underline{x}_{1,d+2})  } \right) (1-\gamma_{ij})\delta_{ij},  s_{k-1}(\underline{x}_{1,d+2},\overline{\gamma}, \overline{\delta})\right\},\nonumber \\
  s_{k}&(&\underline{x}_{1,d+2},\overline{\gamma}, \overline{\delta}) \nonumber\\
    &=& \;\int_{\bbE} \max\Bigg\{\sum_{i,j}\left( 1-p_{ij}^d+ q_{ij}\sum_{m=1}^{d+1}\frac{L_{m}^{i,j}(\underline{x}_{1,d+3}) } { p_{ij}^m L_{0}^{i,j}(\underline{x}_{1,d+3}) } \right)f_{x_{d+2}}^{0,i}(x_{d+3})p_{ij}(1-\gamma_{ij})\delta_{ij}, \nonumber\\
      && \; s_{k-1}(\underline{x}_{2,d+3},\overline{\gamma}, \ovp \circ \ovf0_{x_{d+2}}(x_{d+3})\circ \overline{\delta})\frac{}{} \Bigg\}\mu(dx_{d+3}),\nonumber
\end{eqnarray}
}
where:
{
\setlength\arraycolsep{2pt}
\begin{eqnarray}
s_0(\underline{x}_{1,d+2},\overline{\gamma}, \overline{\delta}) &=& \sum_{i,j}\left( 1-p_{ij}^d + q_{ij}\sum_{m=1}^{d+1}\frac{ L_{m-1}^{i,j}(\underline{x}_{2,d+2}) }{ p_{ij}^m L_{0}^{i,j}(\underline{x}_{2,d+2}) }\right)p_{ij}(1-\gamma_{ij})\delta_{ij}.
\end{eqnarray}
\normalsize
}
Moreover for $k\geq 0$ and vector $\eta_{n+1} = (\underline{X}_{n-d,n+1},\ovPi_{n+1}, \ovBet_{n+1})$, function $s_k$
has the property:
\begin{eqnarray}
\label{PojRozregNieznRozkl-wlasnosc r}
    s_k(\underline{X}_{n-d,n+1},\ovPi_{n+1}, \ovBet_{n+1}) = \frac{ s_k(\underline{X}_{n-d,n+1},\ovPi_{n}, \ovp \circ \ovf0_{X_n}(X_{n+1})\circ \ovBet_{n})}{S(0,\uX_{n,n+1},\ovBet_{n},\ovPi_{n})}.
\end{eqnarray}
\end{lemma}
\begin{pf} Notice that lemmas \ref{PojRozregNieznRozkl-dodatekWzor4}, \ref{PojRozregNieznRozkl-dodatekWzor5},
formulas (\ref{PojRozregNieznRozkl-Pi_n_Rekurencja2}) i (\ref{PojRozregNieznRozkl-B_n_Rekurencja}) enable us to rewrite
function $h(\underline{X}_{n-d,n+1},\ovPi_{n+1}, \ovBet_{n+1})$ in the following way:
{
\small
\setlength\arraycolsep{0pt}
\begin{eqnarray}
\label{PojRozregNieznRozkl-NowaWyplataRekurencja}
h&(&\underline{X}_{n-d,n+1},\ovPi_{n+1}, \ovBet_{n+1}) \\
 &=& \sum_{i,j}\left( 1-p_{ij}^d+ q_{ij}\sum_{m=1}^{d+1}\frac{L_{m}^{i,j}(\underline{X}_{n-d,n+1})} { p_{ij}^m L_{0}^{i,j}(\underline{X}_{n-d,n+1})} \right) (1-\Pi^{i,j}_{n+1})B^{i,j}_{n+1} \nonumber\\
 &=& \sum_{i,j}\left( 1-p_{ij}^d+ q_{ij}\sum_{m=1}^{d+1}\frac{L_{m}^{i,j}(\underline{X}_{n-d,n+1})} { p_{ij}^m L_{0}^{i,j}(\underline{X}_{n-d,n+1})} \right)(1-\Pi^{i,j}(0,\uX_{n,n+1},\Pi_n^{i,j})) \nonumber\\
  && \times \Gamma^{i,j}(0,\uX_{n,n+1},\ovBet_{n},\ovPi_{n}) \nonumber\\
 &=& \sum_{i,j}\left( \frac{(1-p_{ij}^d)p_{ij}(1-\Pi_n^{i,j})B_n^{i,j}}{S(0,\uX_{n,n+1},\ovBet_{n},\ovPi_{n})} f_{X_n}^{0,i}(X_{n+1}) \right. \nonumber\\
\nonumber  & & \left.+ q_{ij}\sum_{m=1}^{d+1}\frac{L_{m-1}^{i,j}(\underline{X}_{n-d,n})} { p_{ij}^m L_{0}^{i,j}(\underline{X}_{n-d,n})}\frac{p_{ij}(1-\Pi_n^{i,j})B_n^{i,j}}{S(0,\uX_{n,n+1},\ovBet_{n},\ovPi_{n})}f_{X_n}^{1,j}(X_{n+1}) \right).
\end{eqnarray}
\normalsize
}
Using definition of operator $\bT$, equation (\ref{PojRozregNieznRozkl-NowaWyplataRekurencja}), for $k=0$ and $(\underline{X}_{n-1-d,n},\ovPi_{n}, \ovBet_{n}) = (\underline{x}_{1,d+2},\overline{\gamma}, \overline{\delta})$ we get
{
\small
\begin{eqnarray}
\label{PojRozregNieznRozkl-wzor na r0}
s_0(\underline{x}_{1,d+2},\overline{\gamma}, \overline{\delta}) &=&\phiE(h(\uX_{n-d,n},X_{n+1},\ovPi_{n+1},\ovBet_{n+1})|\cF_n)\\
\nonumber&\hspace{-10em}=&\hspace{-5em}\int_{\bbE} h(\underline{X}_{n-d,n},y,\ovPi_{n+1}, \ovBet_{n+1})S(0,X_{n},y,\ovBet_{n},\ovPi_{n})\mu(dy) \\
&\hspace{-10em}=&\hspace{-5em}\sum_{i,j} \int_{\bbE} (1-p_{ij}^d)p_{ij}(1-\Pi_n^{i,j}) B_n^{i,j}f_{X_n}^{0,i}(y)\mu(dy) \nonumber\\
\nonumber &\hspace{-10em} &\hspace{-5em} + \sum_{i,j}\int_{\bbE} q_{ij}\sum_{m=1}^{d+1}\frac{L_{m-1}^{i,j}(\underline{X}_{n-d,n})}{ p_{ij}^m L_{0}^{i,j}(\underline{X}_{n-d,n})}p_{ij}(1-\Pi_n^{i,j})B_n^{i,j}f_{X_n}^{1,j}(y)\mu(dy)\\
\nonumber &\hspace{-10em}=&\hspace{-5em}\sum_{i,j} \left( 1-p_{ij}^d + q_{ij}\sum_{m=1}^{d+1}\frac{L_{m-1}^{i,j}(\underline{X}_{n-d,n})}{ p_{ij}^m L_{0}^{i,j}(\underline{X}_{n-d,n})}\right) p_{ij}(1-\Pi_n^{i,j})B_n^{i,j}\nonumber\\
\nonumber &\hspace{-10em}=&\hspace{-5em} \sum_{i,j} \left( 1-p_{ij}^d + q_{ij}\sum_{m=1}^{d+1}\frac{L_{m-1}^{i,j} (\underline{x}_{2,d+2})}{ p_{ij}^m L_{0}^{i,j}(\underline{x}_{2,d+2})}\right)p_{ij}(1-\gamma_{ij})\delta_{ij}.
\end{eqnarray}
\normalsize
}
Hence, applying equations (\ref{PojRozregNieznRozkl-Pi_n_Rekurencja2}) and (\ref{PojRozregNieznRozkl-B_n_Rekurencja}) one more time we end with
{
\setlength\arraycolsep{0pt}
\begin{eqnarray}
\label{PojRozregNieznRozkl-wzor na r01}
    s_0&(&\underline{X}_{n-d,n+1},\ovPi_{n+1}, \ovBet_{n+1}) \nonumber\\
    &=& \sum_{i,j} \left( 1-p_{ij}^d + q_{ij}\sum_{m=1}^{d+1}\frac{L_{m-1}^{i,j}(\underline{X}_{n-d+1,n+1})}{ p_{ij}^m L_{0}^{i,j}(\underline{X}_{n-d+1,n+1})}\right)\frac{p_{ij}(1-\Pi_n^{i,j})p_{ij}f_{X_n}^{0,i}(X_{n+1})B_n^{i,j}}{S(0,\uX_{n,n+1},\ovBet_{n},\ovPi_{n})}\nonumber\\
    &=& \frac{ s_0(\underline{X}_{n-d,n+1},\ovPi_{n}, \ovp \circ \ovf0_{X_n}(X_{n+1})\circ \ovBet_{n})}{S(0,\uX_{n,n+1},\ovBet_{n},\ovPi_{n})}.\nonumber
\end{eqnarray}
}
If $k=1$, then by definition of $\bQ$:
{
\small
\setlength\arraycolsep{0pt}
\begin{eqnarray}
\label{PojRozregNieznRozkl-operator Q}
\!\!\!\!\!\!\!\!  \bQ h(\underline{x}_{1,d+2},\overline{\gamma}, \overline{\delta})&=& \max\Bigl\{\sum_{i,j}\left( 1-p_{ij}^d+ q_{ij}\sum_{m=1}^{d+1}\frac{L_{m}^{i,j}(\underline{x}_{1,d+2})} { p_{ij}^m L_{0}^{i,j}(\underline{x}_{1,d+2}) } \right) (1-\gamma_{ij})\delta_{ij},\\
\!\!\!\!\!\!\!\!  &&\mbox{} s_0(\underline{x}_{1,d+2},\overline{\gamma}, \overline{\delta})\Bigl\}.\nonumber
\end{eqnarray}
\normalsize}
Now, for $(\underline{X}_{n-1-d,n},\ovPi_{n}, \ovBet_{n}) = (\underline{x}_{1,d+2},\overline{\gamma}, \overline{\delta})$, taking into account link between $\Pi_{n-1}^{i,j}$ and $\Pi_{n}^{i,j}$ as well as between $B_{n-1}^{i,j}$ and $B_{n}^{i,j}$ given by (\ref{PojRozregNieznRozkl-Pi_n_Rekurencja2}) and (\ref{PojRozregNieznRozkl-B_n_Rekurencja}), we get with the support of (\ref{expVALUE}):
{
\small
\begin{eqnarray}
\label{PojRozregNieznRozkl-wzor na r1}
s_1(\ux_{1,d+2},\overline{\gamma}, \overline{\delta}) &=& \phiE\Bigl[\max\{h(\underline{X}_{n-d,n},X_{n+1},\ovPi_{n+1}, \ovBet_{n+1}),\\
&&\qquad\qquad s_0(\underline{X}_{n-d,n},X_{n+1},\ovPi_{n+1}, \ovBet_{n+1})\}\mid \cF_n\Bigr] \nonumber\\
&\hspace{-10em}=&\hspace{-5em}\int_{\bbE} \max\left\{\sum_{i,j}\left( 1-p_{ij}^d+ q_{ij}\sum_{m=1}^{d+1}\frac{L_{m}^{i,j}(\underline{X}_{n-d,n},y)} { p_{ij}^m L_{0}^{i,j}(\underline{X}_{n-d,n},y) } \right)\frac{f_{X_n}^{0,i}(y)p_{ij}(1-\Pi_n^{i,j})B_n^{i,j}}{S(0,X_{n},y,\ovBet_{n},\ovPi_{n})},\right. \nonumber\\
&\hspace{-5em}&\hspace{-2em}     \left. \frac{ s_0(\underline{X}_{n-d,n}, y,\ovPi_{n}, \ovp \circ \ovf0_{X_n}(y)\circ \ovBet_{n})}{S(0,X_{n},y,\ovBet_{n},\ovPi_{n})} \right\}S(0,X_{n},y,\ovBet_{n},\ovPi_{n})\mu(dy)\nonumber\\
&\hspace{-10em}=&\hspace{-5em} \int_{\bbE} \max\left\{\sum_{i,j}\left( 1-p_{ij}^d+ q_{ij}\sum_{m=1}^{d+1}\frac{L_{m}^{i,j}(\underline{x}_{2,d+2}, y)} { p_{ij}^m L_{0}^{i,j}(\underline{x}_{2,d+2}, y) } \right)f_{x_{d+2}}^{0,i}(y)p_{ij}(1-\gamma_{ij})\delta_{ij},\right. \nonumber\\
\nonumber &\hspace{-10em}=&\hspace{-5em} \left.  s_0(\underline{x}_{2,d+2},y,\overline{\gamma}, \ovp \circ \ovf0_{x_{d+2}}(y)\circ \overline{\delta})\frac{}{}\right\}\mu(dy).
\end{eqnarray}
\normalsize
}
Basing on (\ref{PojRozregNieznRozkl-wzor na r1}) with the help of (\ref{PojRozregNieznRozkl-Pi_n_Rekurencja2})
and (\ref{PojRozregNieznRozkl-B_n_Rekurencja}) let us verify formula (\ref{PojRozregNieznRozkl-wlasnosc r}):
{
\small
\setlength\arraycolsep{0pt}
\begin{eqnarray}
\label{PojRozregNieznRozkl-wzor na r11}
s_1(\underline{X}_{n-d,n+1},\ovPi_{n+1}, \ovBet_{n+1})&& \nonumber\\
   &\hspace{-20em}=&\hspace{-10em}\int_{\bbE} \max\Bigl\{\sum_{i,j}\left( 1-p_{ij}^d+ q_{ij}\sum_{m=1}^{d+1}\frac{L_{m}^{i,j}(\underline{X}_{n-d+1,n+1},y)} { p_{ij}^m L_{0}^{i,j}(\underline{X}_{n-d+1,n+1},y) } \right) f_{X_{n+1}}^{0,i}(y)p_{ij}(1-\Pi_{n+1}^{i,j})B_{n+1}^{i,j},  \nonumber\\
&\hspace{-20em}&\hspace{-10em} s_0(\underline{X}_{n+1-d,n+1},y,\ovPi_{n+1}, \ovp \circ \ovf0_{X_{n+1}}(y)\circ \ovBet_{n+1})\frac{}{}\Bigr\}\mu(dy)\nonumber\\
&\hspace{-20em}=&\hspace{-10em}\int_{\bbE}\!\! \max\left\{\sum_{i,j}\!\!\left(\!\! 1-p_{ij}^d+ q_{ij}\!\!\sum_{m=1}^{d+1}\frac{L_{m}^{i,j}(\underline{X}_{n-d+1,n+1},y)} { p_{ij}^m L_{0}^{i,j}(\underline{X}_{n-d+1,n+1},y) }\! \right)\! \frac{f_{X_{n+1}}^{0,i}(y)p_{ij}(1-\Pi_{n}^{i,j})B_{n}^{i,j}f_{X_{n}}^{0,i}(X_{n+1})p_{ij}}{S(0,\uX_{n,n+1},\ovBet_{n},\ovPi_{n})}, \right. \nonumber\\
 &\hspace{-20em}&\hspace{-10em} \left.\sum_{i,j}\!\! \left(\!\! 1-p_{ij}^d + q_{ij}\!\!\sum_{m=1}^{d+1}\frac{L_{m}^{i,j}(\underline{X}_{n-d+2,n+1},y)}{ p_{ij}^m L_{0}^{i,j}(\underline{X}_{n-d+2,n+1},y)}\!\!\right)\!\!\frac{p_{ij}(1-\Pi_n^{i,j})p_{ij}f_{X_n}^{0,i}(X_{n+1})p_{ij}f_{X_{n+1}}^{0,i}(y)B_n^{i,j}}{S(0,\uX_{n,n+1},\ovBet_{n},\ovPi_{n})}\!\!\right\}\mu(dy)\nonumber\\
    &=& \frac{ s_1(\underline{X}_{n-d,n+1},\ovPi_{n}, \ovp \circ \ovf0_{X_n}(X_{n+1})\circ \ovBet_{n})}{S(0,\uX_{n,n+1},\ovBet_{n},\ovPi_{n})}.\nonumber
\end{eqnarray}
\normalsize
}
Suppose that lemma \ref{PojRozregNieznRozkl-OperatoryDlaNowejWyplaty} holds for some $k >1$. We will show that equations characterizing $\bQ^{k+1}h$ and $s_{k+1}$ are true and that condition (\ref{PojRozregNieznRozkl-wlasnosc r}) for $s_{k+1}$ is satisfied. It follows from definition of operator $\bQ^{k+1}$ that:
\begin{eqnarray}
\label{PojRozregNieznRozkl-operator Qk}
  \bQ^{k+1} h&(&\underline{x}_{1,d+2},\overline{\gamma}, \overline{\delta}) \\
  &=& \max\Big\{\sum_{i,j}\left( 1-p_{ij}^d+ q_{ij}\sum_{m=1}^{d+1}\frac{L_{m}^{i,j}(\underline{x}_{1,d+2})} { p_{ij}^m L_{0}^{i,j}(\underline{x}_{1,d+2}) } \right) (1-\gamma_{ij})\delta_{ij},\nonumber\\ &&s_k(\underline{x}_{1,d+2},\overline{\gamma}, \overline{\delta})\Big\}. \nonumber
\end{eqnarray}
Given $(\underline{X}_{n-1-d,n},\ovPi_{n}, \ovBet_{n}) = (\underline{x}_{1,d+2},\overline{\gamma}, \overline{\delta})$ and basing on inductive assumption we have also:
{
\small
\begin{eqnarray}
\label{PojRozregNieznRozkl-wzor na rk+1}
     s_{k+1}(\underline{x}_{1,d+2},\overline{\gamma}, \overline{\delta}) &=& \phiE\Bigl[\max \Bigl\{h(\underline{X}_{n-d,n},X_{n+1},\ovPi_{n+1}, \ovBet_{n+1}),\\
     &&\mbox{} s_k(\underline{X}_{n-d,n},X_{n+1},\ovPi_{n+1}, \ovBet_{n+1})\Bigr\} \mid \cF_n\Bigr]\nonumber \\
    &\hspace{-10em}=&\hspace{-5em} \int_{\bbE} \max\left\{\sum_{i,j}\left( 1-p_{ij}^d+ q_{ij}\sum_{m=1}^{d+1}\frac{L_{m}^{i,j}(\underline{X}_{n-d,n},y)} { p_{ij}^m L_{0}^{i,j}(\underline{X}_{n-d,n},y) } \right)\frac{f_{X_n}^{0,i}(y)p_{ij}(1-\Pi_n^{i,j})B_n^{i,j}}{S(0,X_{n},y,\ovBet_{n},\ovPi_{n})},\right. \nonumber\\
      && \left. \frac{ s_k(\underline{X}_{n-d,n}, y,\ovPi_{n}, \ovp \circ \ovf0_{X_n}(y)\circ \ovBet_{n})}{S(0,X_{n},y,\ovBet_{n},\ovPi_{n})} \right\}S(0,X_{n},y,\ovBet_{n},\ovPi_{n})\mu(dy)\nonumber\\
    &\hspace{-10em}=&\hspace{-5em}\int_{\bbE} \max\left\{\sum_{i,j}\left( 1-p_{ij}^d+ q_{ij}\sum_{m=1}^{d+1}\frac{L_{m}^{i,j}(\underline{X}_{n-d,n}, y)} { p_{ij}^m L_{0}^{i,j}(\underline{X}_{n-d,n}, y) } \right)f_{X_n}^{0,i}(y)p_{ij}(1-\Pi_n^{i,j})B_n^{i,j},\right. \nonumber\\
      && \left.  s_k(\underline{X}_{n-d,n},y,\ovPi_{n}, \ovp \circ \ovf0_{X_n}(y)\circ \ovBet_{n})\frac{}{}\right\}\mu(dy)\nonumber\\
    &\hspace{-10em}=&\hspace{-5em} \int_{\bbE} \max\left\{\sum_{i,j}\left( 1-p_{ij}^d+ q_{ij}\sum_{m=1}^{d+1}\frac{L_{m}^{i,j}(\underline{x}_{2,d+2}, y)} { p_{ij}^m L_{0}^{i,j}(\underline{x}_{2,d+2}, y) } \right)f_{x_{d+2}}^{0,i}(y)p_{ij}(1-\gamma_{ij})\delta_{ij},\right. \nonumber\\
\nonumber      && \left.  s_k(\underline{x}_{2,d+2},y,\overline{\gamma}, \ovp \circ \ovf0_{x_{d+2}}(y)\circ \overline{\delta})\frac{}{}\right\}\mu(dy).
\end{eqnarray}
\normalsize
}
Finally, using(\ref{PojRozregNieznRozkl-wzor na rk+1}) we obtain:
{
\small
\setlength\arraycolsep{0pt}
\begin{eqnarray}
\label{PojRozregNieznRozkl-wzor na rk1}
    s_{k+1}&(&\underline{X}_{n-d,n+1},\ovPi_{n+1}, \ovBet_{n+1}) \nonumber\\
    &=& \int_{\bbE} \max\left\{\sum_{i,j}\left( 1-p_{ij}^d+ q_{ij}\sum_{m=1}^{d+1}\frac{L_{m}^{i,j}(\underline{X}_{n-d+1,n+1},y)} { p_{ij}^m L_{0}^{i,j}(\underline{X}_{n-d+1,n+1}, y) } \right) f_{X_{n+1}}^{0,i}(y)p_{ij}(1-\Pi_{n+1}^{i,j})B_{n+1}^{i,j}, \right. \nonumber\\
      && \left. s_{k}(\underline{X}_{n+1-d,n+1},y,\ovPi_{n+1}, \ovp \circ \ovf0_{X_{n+1}}(y)\circ \ovBet_{n+1})\frac{}{}\right\}\mu(dy)\nonumber\\
    &=& \int_{\bbE}\!\! \max\left\{\sum_{i,j}\!\!\left(\!\! 1-p_{ij}^d+ q_{ij}\!\!\sum_{m=1}^{d+1}\frac{L_{m}^{i,j}(\underline{X}_{n-d+1,n+1}, y)} { p_{ij}^m L_{0}^{i,j}(\underline{X}_{n-d+1,n+1},y) }\! \right)\! \frac{f_{X_{n+1}}^{0,i}(y)p_{ij}(1-\Pi_{n}^{i,j})B_{n}^{i,j}f_{X_{n}}^{0,i}(X_{n+1})p_{ij}}{S(0,\uX_{n,n+1},\ovBet_{n},\ovPi_{n})}, \right. \nonumber\\
      && \left. \frac{s_{k}(\underline{X}_{n+1-d,n+1},y,\ovPi_{n}, \ovp \circ \ovf0_{X_{n+1}}(y)\circ  \ovp \circ \ovf0_{X_{n}}(X_{n+1})\circ\ovBet_{n})}{S(0,\uX_{n,n+1},\ovBet_{n},\ovPi_{n})}\right\}\mu(dy)\nonumber\\
    &=& \frac{ s_{k+1}(\underline{X}_{n-d,n+1},\ovPi_{n}, \ovp \circ \ovf0_{X_n}(X_{n+1})\circ \ovBet_{n})}{S(0,\uX_{n,n+1},\ovBet_{n},\ovPi_{n})},\nonumber
\end{eqnarray}}
\end{pf}

%\bibliographystyle{plainnat}%elsart-harv}%gSSR}
%\bibliography{DisWSKSz28x09}%%gSSRguide
%\markboth{Sarnowski \& Szajowski}{Sequential Analysis}
\end{document}